\documentclass{article}
\usepackage[T2A]{fontenc}
\usepackage[cp1251]{inputenc}
\usepackage[english,russian]{babel}
\usepackage[tbtags]{amsmath}
\usepackage{amsfonts,amssymb,mathrsfs,amscd}

\let\No=\textnumero

\usepackage[hyper]{izv2e}
\JournalName{СЕРИЯ МАТЕМАТИЧЕСКАЯ}

\numberwithin{equation}{section}
\theoremstyle{plain}
\newtheorem{theorem}{Теорема}
\newtheorem{lemma}{Лемма}[section]

\theoremstyle{definition}

\def\klabel#1{}

\begin{document}
	 
\udk{517.984.50}

\date{}

\author{Р.\,С.~Сакс}
\address{Институт математики c ВЦ УФИЦ РАН, \\ ул.Чернышевского, 112, 450077, г. Уфа, Россия
}
\email{romen-saks@yandex.ru}

\title{Операторы вихрь и градиент дивергенции в пространствах Соболева}

\markboth{Р.\,С.~Сакс}{Операторы $\mathrm{ROT}$ и $\nabla\,\mathrm{DIV}$  	в пространствах Соболева}

\maketitle

\begin{fulltext}

\begin{abstract} С.Л. Соболев, 	изучая  краевые задачи для  полигармонического уравнения $\Delta^m\,u=\rho$ в пространствах ${W} _2^{(m)}(\Omega)$ c обобщённой правой частью, заложил фундамент теории этих  пространств \cite{sob}  $\S 9$ гл. 12.
		
	Операторы   градиент дивергенции  и  ротор ротора  ($\nabla \text{div}$   и	 	$ \text{rot}^2$)  и их степени являются аналогами скалярного оператора 	$\Delta^m$  в   ортогональных подпространствах  $\mathcal{A}$  и $\mathcal{B}$  в  $\mathbf{L}_{2}(G)$.
	Они порождают разновидности $ \mathbf{A}^{2k}(G)$ и 
	 $\mathbf{W}^m(G)$ пространств Соболева  потенциальных и вихревых полей,  а их прямые суммы  $ \mathbf{A}^{2k}(G) \oplus \mathbf{W}^m(G)$ -  сеть пространств.  Её элементы,	 
	  классы $ \mathbf{C}(2k, m)\equiv \mathbf{A}^{2k}\oplus \mathbf{W}^m$, играют роль 	${W} _2^{(m)}(G)$ 	в ограниченной области   $ G$  в 
	 $\mathbb{R}^3$ с гладкой границей $\Gamma$. 
	 
	 В этой статье мы рассмотрим  свойства операторов  градиент дивергенции  и  ротор и построим пространства 
	 $ \mathbf{A}^{2k}(G)$ и $\mathbf{W}^m(G)$. 
	 
	 Приведём работы физиков о появлении безсиловых  полей Бельтрами в астрофизике  и в физике плазмы.
	 
	 Краевые задачи  для операторов  $\nabla \text{div}+\lambda I$ 	 	   и	 	$ \text{rot}+\lambda I$ в  классах $  \mathbf{C}(2k, m)$  рассмотрим в следующей работе.
	 
	 \end{abstract} 
	
Подобно тому, как течения жидкости разделяют на ламинарные и турбулентные, векторные поля в 	$\mathbf{L}_{2}(G) $  разделяются на потенциальные (безвихревые) и соленоидальные (вихревые).
	
		Пространство 	$\mathbf{L}_{2}(G) $  разлагается на ортогональные подпространства 	
 $\mathcal{A}$  и $\mathcal{B}$: 	$\mathbf{L}_{2}(G)=\mathcal{A}\oplus \mathcal{B}$. 	
	 В свою очередь   
	$ \mathcal{A}= \mathcal{A}_H\oplus \mathbf{A}^0$ и $ \mathcal{B}=\mathcal{B}_H \oplus \mathbf{V}^0$,
	где 	$\mathcal{A}_H $ и $\mathcal{B}_H $- нуль-пространства операторов 	$\nabla \text{div}$   в $\mathcal{A}$  и  $ \text{rot}$  в  $\mathcal{B}$; 	
	размерности   $\mathcal{A}_H $ и $\mathcal{B}_H $ конечны и определяется топологией границы. 

	Собственные поля	
	 оператора	$\nabla \text{div}$  (соотв., 	$ \text{rot}$)  с ненулевыми  собственными значениями  
	используются при построении  ортонормированного базиса  в  $ \mathbf{A}^0$  (соотв.,   в $\mathbf{V}^0$). 	 
	
		 Операторы   $\nabla\mathrm{div}$ и  $\mathrm{rot}$ аннулируют друг друга	  и 	проектируют $\mathbf{L}_{2}(G) $   на   $\mathcal{A}$  и $\mathcal{B}$,  причем $\mathrm{rot}\, \mathbf {u}=0$ при $\mathbf{u}\in	\mathcal{A}$,  а $\nabla\mathrm div \mathbf {v}=0$ при $\mathbf{v}\in	\mathcal{B}$  \cite{hw}.	
		 	
	 Лапласиан в  $\mathbf{L}_{2}(G)$
	выражается через них: 
	$\mathrm{\Delta} \mathbf {v}
	\equiv \nabla \mathrm{div}\,\mathbf {v}
	-(\mathrm{rot})^2\, \mathbf {v}$.	
	Поэтому 
	$\Delta^m\equiv  (\nabla \mathrm{div})^m$   в $   \mathcal{A}$  и  $ \Delta^m  \equiv(-1)^m(\mathrm{rot})^{2m}$    в       $ \mathcal{B}$  при   $ m\geq  1$.

	Аналоги 	 пространств  	${W}_2^{(m)}(G)$ в  классах  $\mathcal{A}$  и $\mathcal{B}$  -	это пространства   $\mathbf{A}^{2k}(G)$ и 	$\mathbf{W}^m(G)$,   порядков   $2k>0$  и $ m>0$,  а      $\mathbf{A}^{-2k}(G)$ и 	$\mathbf{W}^{-m}(G)$- двойственные им  пространства,  сопряжённые с  $\mathbf{A}_0^{2k}(G)$ и 	$\mathbf{W}_0^m(G)$.
	 По определению    
	\[ \mathbf{A}^{2k}(G)=
	\{\mathbf{f}\in\mathbf{A}^0_{\gamma}(G),...,
	(\nabla\mathrm{div})^k\,\mathbf{f}\in\mathbf{A}^0_{\gamma}(G) \}\quad
\quad \text{при}\quad \, k\geq 1,    \eqno{(0.1)}\]

	\[\mathbf{W}^{m}(G)=
	\{\mathbf{f}\in\mathbf{V}^{0}(G),..., \mathrm{rot}^m
	\mathbf{f}\in\mathbf{V}^{0}(G) \}\quad \text{при}\quad m\geq
	1.   \eqno{(0.2)}\]
	
Они образуют две шкалы (цепи) вложенных пространств:
\[\subset \mathbf{A}^{2k}\subset...\subset\mathbf{A}^{2}\subset  \mathbf{A}^0\subset \mathbf{A}^{-2}\subset...\subset \mathbf{A}^{-2k}\subset \eqno{(0.3)}\]
\[\subset \mathbf{W}^{m}\subset...\subset \mathbf{W}^1\subset \mathbf{V}^0\subset  \mathbf{W}^{-1}\subset...\subset \mathbf{W}^{-m}\subset     \eqno{(0.4)}\]

В них   действуют  операторы  $\mathcal{N}_d$ и $S$ - самосопряженные расширения   операторов	$\nabla \text{div}$ и  $ \text{rot}$  в   пространства  $ \mathbf{A}^0$ и  $\mathbf{V}^{0}$  безвихревых и  вихревых полей.

    $\mathcal{N}_d$ и $S$ отображают пространства   $ \mathbf{A}^{2k}$ на  $ \mathbf{A}^{2(k-1)}$  и  $ \mathbf{W}^{m}$ на  $ \mathbf{W}^{m-1}$, соответственно,  а   операторы    $\mathcal{N}^{-1}_d$  и $S^{-1}$ -- в обратную сторону.

Отображения 
$\mathcal{N}_d^{2k}:\mathbf{A}^{2k}(G)\to \mathbf{A}^{-2k}(G)$, \quad 
$S^{2m}: \mathbf{W}^m\to \mathbf{W}^{-m}$ и  обратные  отображения 
 $\mathcal{N}_d^{-2k}$  при  $k\geq 1$  и
$S^{-2m}$   при $m\geq 1$ также рассмотрены.  
Доказано, что
{\it 	
	уравнение $(\nabla \mathrm{div})^{2k}\,\mathbf{u}=\mathbf{v}$ 
	при заданном $\mathbf{v}$  в 	объединении  $ \mathbf{A}^{-2n}$ и $k\geq 1$ разрешимо в
	пространстве $\mathbf{A}^{2k}$ тогда и только тогда, когда
	$\mathbf{v}\in\mathbf{A}^{-2k}$.\quad
	Его решение   $\mathbf{u}= \mathcal{N}_d^{-2k}\mathbf{v}$
	в  фактор-пространстве	$\mathcal{A}/\mathcal{A}_H$   определяется однозначно.}

Аналогично, {\it 	при заданном $\mathbf{v}$  в 
	объединении  $\mathbf{W}^{-n}$
	уравнение $\mathrm{rot}^{2m}\,\mathbf{u}=\mathbf{v}$  разрешимо в
	пространстве $\mathbf{W}^m(G)$ тогда и только тогда, когда
	$\mathbf{v}\in \mathbf{W}^{-m}(G)$.\,
	Его решение $\mathbf{u}=S^{-2m}\mathbf{v}$ в  классе смежности   
	$ \mathcal{B}/\mathcal{B}_H$  определяется 
	однозначно. }

Пары пространств из цепочек (0.3) и (0.4) образуют сеть пространств Соболева,  её элементы - классы  $ \mathbf{C}(2k, m)(G)\equiv \mathbf{A}^{2k}(G) \oplus \mathbf{W}^m(G)$;
 класс  $  \mathbf{C}(2k, 2k)$ совпадает с пространством Соболева    $\mathbf{H}^{2k}(G)$. 
Они принадлежат     $\mathbf{L}_{2}(G)$, если  $k\geq 0$ и
$m\geq 0$.

	Открылось широкое поле задач: изучение операторов  $(\mathrm{rot})^p$, $ (\nabla \, \mathrm{div})^p$  при  $p=1,2, ...$ и других в классах   Соболева $ \mathbf{C}(2k, m)$. 
	
	Мы рассмотрим их  в следующих  работах.

 В работе использованы обозначения и результаты работ автора
  \cite{saVS20, ds18,sa15a,sa15b}.  
Библиография: 33 названия.

*
 
 

\section { Основные подпространства $\mathbf{L}_{2}(G)$}

 Мы рассматриваем линейные
пространства над полем $\mathbb{R}$ действительных чисел. Через
$\mathbf{L}_{2}(G)$ обозначаем пространство Лебега вектор-функций  (полей), квадратично интегрируемых в $G$ с внутренним произведением

 $(\mathbf {u},\mathbf {v})= \int_G \mathbf
{u}\cdot \mathbf {v}\,d \mathbf {x}$ и нормой
$\|\mathbf{u}\|= (\mathbf {u},\mathbf {u})^{1/2}$.
\subsection {Шкала  пространств Соболева}
 Пространство Соболева, состоящее из полей,
 принадлежащих $\mathbf{L}_{2}(G)$ вместе с обобщенными производными
 до порядка $ m> 0$, обозначается через
$\mathbf{H}^{m}(G)$, $\|\mathbf {f}\|_m$ -норма его элемента
$\mathbf {f}$;

$\mathbf{H}^{0}(G)\equiv\mathbf{L}_{2}(G)$. 
 $\mathbf{H}^m(G)$ - гильбертово пространство со 
  скалярным произведением:
\[	 (\mathbf{f},\mathbf{g})_m=(\mathbf{f},\mathbf{g})+
\int_{G} \sum_{|\alpha|=m}\frac{m!}{\alpha
	!}\partial^{\alpha}\mathbf{f}\cdot\partial^{\alpha}\mathbf{g}
d \mathbf{x},  \quad   \|\mathbf {f}\|_m^2= (\mathbf{f},\mathbf{f})_m.     \eqno{(1.1)}\] 	

Замыкание в норме $\mathbf{H}^{m}(G)$ множества $[\mathcal{C}^{\infty}_0(G)]^3$
обозначается через $\mathbf{H}^{m}_0(G)$.

 Двойственное пространство Соболева
отрицательного порядка $\mathbf{H}^{-m}(G)$ сопряжено с
$\mathbf{H}^{m}_0(G)$

 На лекции    в  НГУ  в 1962  году  С.Л.Соболев рисовал всю цепь  вложенных пространств: 
 \[\subset \mathbf{H}^{m}\subset...\subset \mathbf{H}^1\subset \mathbf{L}_2\subset  \mathbf{H}^{-1}\subset...\subset \mathbf{H}^{-m} \subset       \eqno{(1.2)}\]
Он обозначал их  $W_2^{(m)}(G)$ в  \cite{sob}  $\S 9$ гл. 12. 

 Их обозначают также $\mathbf{H}^{m}(G)$ (см.  В.П.Михайлов \cite{mi} $\S 4$ гл. 3).
 
 В области $G$ с гладкой границей $\Gamma$ в каждой точке $y\in\Gamma$
 определена нормаль $\mathbf {n}(y)$ к $\Gamma$.
Поле $\mathbf {u}$ из $\mathbf{H}^{m+1}(G)$ имеет след $
\gamma(\mathbf {n}\cdot\mathbf {u})$ на $\Gamma$ его нормальной
компоненты, который принадлежит пространству Соболева-Слободецкого
$\mathbf{H}^{m+1/2}(G)$, $|\gamma(\mathbf {n}\cdot\mathbf
{u})|_{m+1/2}$- его норма.

\subsection {Пространства потенциальных и соленоидальных
	$\mathcal{A}$ и $\mathcal{B}$ в $\mathbf{L}_{2}(G)$}
Пусть $h$- функция  из $H^{1}(G)$, а $\mathbf{u}=\nabla h$ - ее
градиент. 
По определению ${\mathcal{{A}}}(G) =\{\nabla h, h\in H^1(G)\}$,  
 а  $\mathcal{B}$ - ортогональное дополнение  $\mathcal{A}$ в пространстве  $\mathbf{L}_{2}(G)$.  

 Из соотношений ортогональности  $(\mathbf {u},\nabla h)=0$ для любой   $ h\in H^1(G)$ при $\mathbf{u}\in\mathbf{H}^{1}(G)$ вытекает, что $ \mathrm{div} \mathbf{u}=0$  в  $G$, \,$\gamma(\mathbf{n}\cdot \mathbf{u})=0 $.  Поэтому  ${\mathcal{{B}}}(G)$ обозначают ещё так:  $\mathcal{B}(G)=\{\mathbf{u}\in\mathbf{L}_{2} (G): \mathrm{div} \mathbf{u}=0 \quad \text{в}\quad G, \,\gamma(\mathbf{n}\cdot \mathbf{u})=0 \}$.  
 \footnote{Если   $\mathbf{u}$  и  $ \mathrm{div}\mathbf{u}\in\mathbf{L}_{2} (G)$,  то след $\,\gamma(\mathbf{n}\cdot \mathbf{u})$ существует \cite{rt}.} 
 Итак, 
\[\mathbf{L}_{2}(G)=
{\mathcal{{A}}}(G)\oplus{\mathcal{{B}}}(G).  \eqno{(1.3)}\]
{\it Замечание.} {\small Это разложение содержится в статье   Z.Yoshida и Y.Giga    \cite{yogi}}.   Авторы  называют его разложением Вейля \cite{hw}, а  ${\mathcal{{B}}}(G)$  обозначают как ${L}_{\sigma}^2(G)$.

 Если граница $\Gamma$ имеет положительный род $\rho$,  то 
   $ \mathcal{A}$  содержит  подпространство 
   \[\mathcal{A}_H=\{\mathbf{v
	}\in\mathbf{L}_{2}(G):\,\nabla\text{div}\mathbf{v}=0, \quad
 \text{rot}\,\mathbf{v}=0 \quad \text{в}\quad G, \quad \gamma (\mathbf{n}\cdot \mathbf{v})=0 \},   \eqno{(1.4)}\]  
а  $ \mathcal{B}$ содержит  подпространство безвихревых соленоидальных полей
\begin{equation*}\klabel{bh 1} \mathcal{B}_H=\{\mathbf{u}\in\mathbf{L}_{2}(G):\,\mathrm{div} \mathbf{u}=0, \,\, \mathrm{rot} \mathbf{u}=0  \quad \text{в}\quad G, \quad \gamma(\mathbf{n}\cdot \mathbf{u})=0 \}.  \eqno(1.5) \end{equation*}
  Размерность $\mathcal{B}_H$ равна $\rho$ \cite{boso} и
его     базисные поля
 $\mathbf{h}_j\in \mathbf{C}^{\infty}\,(\bar{G})$,  $j=1,..,\rho$ \cite{hw}. 
 Размерность $\mathcal{A}_H$  не меньше  $\rho$, так как 
 $\mathcal{B}_H\subset \mathcal{A}_H$.  Его     базисные поля
$\mathbf{g}_l \in \mathbf{C}^{\infty}\,(\bar{G})$,  $l=1,..,\rho_1\geq\rho$ (см. п. 1.7).

Отметим, что    у сферы
 размерность $\rho =\dim\mathcal{B}_H =0$,  у тора $\rho=1$  и
 $\rho\geq 1$  у сферы с ручками (числом  $\rho$)  . 

 Ортогональное дополнение   в $\mathcal{A}$ к   $\mathcal{A}_H$  
 обозначается    $\mathbf {A}^{0} (G)$.
 
 Ортогональное дополнение   в $\mathcal{B}$  к   $\mathcal{B}_H$ обозначается  $\mathbf {V}^{0} (G)$   и называется классом {\it вихревых}\,   полей   \cite{saVSTU}. Так что   
\begin{equation*}\klabel{bhr 1}
\mathcal{A}(G)=\mathcal{A}_{H} (G) \oplus \mathbf {A} ^{0} (G), \quad 
\mathcal{B}(G)=\mathcal{B}_{H} (G) \oplus \mathbf {V} ^{0} (G). 
 \eqno(1.6) \end{equation*}
{\it  В шаре   $B$,   множества $\mathcal{A}_H$ и $\mathcal{B}_H$   пусты  и  $\mathbf {A}^{0} (B)= \mathcal{A}(B)$,  а   $\mathbf {V}^{0} (B)=\mathcal{B}(B)$.}

{\it Замечание. }  
 С.Л.Соболев\,\cite{sob54}, \,О.А.Ладыженская\,\cite{lad}, К.Фридрихс \,\cite{fri},  Э.Б.Быховский и Н.В. Смирнов \, \cite{bs}  также	приводят  разложения $\mathbf{L}_{2}(G)$ на   	ортогональные подпространства. 
Причём,  С.Л.~Соболев  	предполагает, что область  (он  обозначает её $\Omega$)  гомеоморфна шару.    Z.Yoshida и Y.Giga отмечают в  \cite{yogi},  что разложение $ \mathcal{B}(G)$  (1.6) содержится в книге C.B. Morrey \cite{mo}.

  Мы   будем придерживаться разложений (1.3), (1.6). 
  
    \subsection{Операторы   $\nabla\mathrm{div}$   и   $\mathrm{rot}$ - проекторы $\mathbf{L}_{2}(G)$ на 	$\mathcal{A}$ и $\mathcal{B}$} Операторы градиент,  ротор (вихрь) и дивергенция определяются в трехмерном векторном анализе, например,  в курсе В.А.Зорича  \cite{zo}.      Им соответствует оператор $d$ внешнего
дифференцирования на формах $\omega^k$ степени $k=0,1$ и 2.
Соотношения $dd\omega^k=0$ при $k=0,1$ имеют вид $\mathrm{rot}\,\nabla h=0$ и $\mathrm{div}\, \mathrm{rot} \mathbf{u}=0$ для гладких функций $h$ и  $ \mathbf{u}$.
Следовательно,   операторы   $\nabla\mathrm{div}$ и  $\mathrm{rot}$ аннулируют друг друга: 
\[\nabla\mathrm{div}\, \mathrm{rot}\, \mathbf{u}=0, \quad     \mathrm{rot}\,\nabla\mathrm{div} \mathbf{u}=0.     \eqno{(1.7)}\]    
  Оператор Лапласа
выражается через них  и скалярный оператор $\Delta_c$: 
\begin{equation*}\klabel{dd 1}\mathrm{\Delta} \mathbf {v}
\equiv \nabla \mathrm{div}\,\mathbf {v}
-(\mathrm{rot})^2\, \mathbf {v}= \Delta_c I_3\,\mathbf {v} , \quad  
 \,\,  \mathbf{v}=(v_1, v_2, v_3),  \quad  \Delta_c v_j\equiv  \mathrm{div} \nabla {v_j}\eqno(1.8) \end{equation*}
 Оператор Лапласа эллиптичен  \cite{vol, so71, sa75},  а операторы $\mathrm{rot}$ и $\nabla \mathrm{div}$ не являются эллиптическими.  Они
вырождены, причем $\mathrm{rot}\, \mathbf {u}=0$ при $\mathbf{u}\in
\mathcal{A}$,  а $\nabla\mathrm div \mathbf {v}=0$ при $\mathbf{v}\in
\mathcal{B}$ в смысле  $\mathbf{L}_{2}(G)$ \cite{hw}.   
Поэтому 
\[\Delta \mathbf{v}\equiv  \nabla \mathrm{div} \mathbf{v}   \quad  \text{при}    \quad  \mathbf{v}\in \mathcal{A},  \quad  \Delta \mathbf{u} \equiv -\mathrm{rot}\, \mathrm{rot}\mathbf{u}  \quad  \text{при}  \quad      \mathbf{u}\in \mathcal{B}.   \eqno(1.9)\]

\subsection{Краевые   задачи для $\text{rot}$    и  $ \nabla\text{div}$ в  пространствах Соболева} 
  В классе Б.Вайнберга и В.Грушина  равномерно неэллиптических псевдодифференциальных операторов  \cite{vagr}, автор выделил в  \cite{saVS20}  подкласс [REES p] обобщённо эллиптических дифференциальных операторов   и доказал, что  {\it операторы $ \text{rot}+\lambda \, I$    и  $ \nabla\text{div}+\lambda \,I$   первого и второго порядков  при $\lambda\neq 0$ принадлежат  классу} [REES 1]. В пространствах Соболева $\mathbf{H}^{s}(G)$  изучены     краевые   задачи.  Им  соответствуют операторы $\mathbb{A}$ и $\mathbb{B}$, которые расширятся до эллиптических по В.Солонникову    переопределённых операторов $\mathbb{A}_R$ и $\mathbb{B}_R$,  ограниченных в  пространствах $\mathbf{H}^{s}(G)$ при целом $s\geq 0$:
\begin{equation*}\klabel{op 1r} \mathbb{A}_R\mathbf{u}\equiv\left( \begin{matrix}
\mathrm{rot} +\lambda I \\
\lambda \, \mathrm{div}\\
\gamma
\mathbf{n}\cdot \end{matrix}\right)\mathbf{u}:
\mathbf{H}^{s+1}(G)\rightarrow\left(
\begin{matrix}\mathbf{H}^{s}(G)\\
H^s(G)\\
H^{s+1/2}(\Gamma)\end{matrix}\right),    \eqno(1.10)
\end{equation*}
\begin{equation*}\klabel{op 1r} \mathbb{B}_R\mathbf{u}\equiv\left( \begin{matrix}
\nabla\,\mathrm{div} +\lambda I \\
\lambda \, \mathrm{rot}\\
\gamma
\mathbf{n}\cdot \end{matrix}\right)\mathbf{u}:
\mathbf{H}^{s+2}(G)\rightarrow\left(
\begin{matrix}\mathbf{H}^{s}(G)\\
\mathbf{H}^{s+1}(G)\\
H^{s+3/2}(\Gamma)\end{matrix}\right).               \eqno(1.11)
\end{equation*}

Из  Теоремы 1.1  В.Солонникова \cite{so71}  в работе \cite{saVS20} получена:
 \begin{theorem} \klabel{rot   1}
	Оператор $\mathbb{A}_R$	имеет левый регуляризатор.  
	Его ядро конечномерно и для любых $ \mathbf{u}\in  \mathbf{H}^{s+1}(G)$ и $ \lambda \neq 0 $ ( с  постоянной  $C_s =C_s(\lambda )>0$,  зависящей только от $s, \lambda $)	выполняется оценка: 
		\begin{equation*}
	\klabel{arot__1_} C_s\|\mathbf{u}\|_{s+1}
	\leq\|\mathrm{rot} \mathbf{u}\|_{s}+
	|\lambda|\|\mathrm{div} \mathbf{u}\|_{s}+
	|\gamma({\mathbf{n}}\cdot\mathbf{u})|_{s+1/2}+
	\|\mathbf{u}\|_{s}.\eqno{(1.12)}
	\end{equation*}
\end{theorem} 
\begin{theorem} \klabel{Nd_1}
	Оператор $\mathbb{B}_R$ 	имеет левый регуляризатор.
	Его ядро конечномерно и   для любых $\mathbf{v}\in  \mathbf{H}^{s+2}(G)$ и $\lambda \neq 0 $ ( с  постоянной  $C_s =C_s(\lambda )>0$,  зависящей только от $s, \lambda $)	выполняется оценка:
	\begin{equation*}
	\klabel{ond__2_} C_s\|\mathbf{v}\|_{s+2}
	\leq|\lambda|\|\mathrm{rot} \mathbf{v}\|_{s+1}+
	\|\nabla\mathrm{div} \mathbf{v}\|_{s}+
	|\gamma({\mathbf{n}}\cdot\mathbf{v})|_{s+3/2}+ \|\mathbf{v}\|_{s}.\eqno{(1.13)}
	\end{equation*}
\end{theorem}

 Топологических ограничений на область нет, предполагается ее
связность, ограниченность и гладкость границы.
Оценка (1.12) известна давно (см. \cite{fri, yogi}). 
 Здесь мы показываем, что  для  операторов класса  [REES p]  аналогичные оценки
легко получать из Теоремы В.Солонникова. 

 Формулы $\mathbf{u}\cdot\nabla h+ h\mathrm{div}\mathbf{u}=\mathrm{div}(h \mathbf{u}) $, \quad 
 $\mathbf{u}\cdot\mathrm{rot} \mathbf{v}- \mathrm{rot}
 \mathbf{u}\cdot\mathbf{v}=\mathrm{div}[\mathbf{v},\mathbf{u}]$, где
 $\mathbf{u}\cdot\mathbf{v}$ и 
 $[\mathbf{v},\mathbf{u}]$ -скалярное и  векторное произведения в   $R^3$, и интегрирование по  $G$ используются при определении операторов $\nabla\,\mathrm{div}$ и $\mathrm{rot}$ в $\mathbf{L}_{2}(G)$.  
  Интегрируя и применяя формулу Гаусса-Остроградского, имеем
 \[ \int_G [\mathrm{rot} \mathbf{u}\cdot \mathbf{v}
 -\mathbf{u}\cdot \mathrm{rot} \mathbf{v}]d\mathbf{x} =  \int_{\Gamma}   \mathbf{n}\cdot [\mathbf{v},\mathbf{u}]d\mathbf{S}.  \eqno{(1.14)}\] 
 \[ \int_G [\nabla \mathrm{div} \mathbf{u}\cdot \mathbf{v}
 -\mathbf{u}\cdot \nabla\mathrm{div} \mathbf{v} ]d\mathbf{x} =  \int_{\Gamma} [ (\mathbf{n}\cdot \mathbf{v})\mathrm{div} \mathbf{u}+  (\mathbf{n}\cdot \mathbf{u})\mathrm{div} \mathbf{v}]d\mathbf{S}.  \eqno{(1.15)}\] 
  \subsection{Операторы  $S$ и  $ \mathcal{N}_d$ -  самосопряженные расширения     $\text{rot}$ и  $\nabla\text{div}$ в  $\mathbf{L}_{2}(G)$} 
  Пусть  $ \mathcal{A}_{\gamma} (G) = \{\nabla\,h, h\in H^2(G): \gamma (\mathbf{n}\cdot \nabla) h =0 \}$,  $\mathbf{A}^0_{\gamma} =\mathbf{A}^0\cap \mathcal{A}_{\gamma}$.
  
  Области определения операторов  $S$ и  $ \mathcal{N}_d$ - это пространства:
 \[   \mathbf{W}^1= \{\textbf{ f}\in  \mathbf{V}^0, \,\,  \text{rot}\textbf{ f}\in  \mathbf{V}^0 \} \quad   \text{и} \quad        \mathbf{A}^{2}=\{\textbf{ f}\in  \mathbf{A}^0_{\gamma}, \,\, \nabla  \text{div}\textbf{ f}\in  \mathbf{A}^0_{\gamma} \}, \eqno{(1.16)}\]	  	
и 
$S\mathbf{u}=\text{rot}\mathbf{u}$ при  $\mathbf{u} \in \mathcal{D}(S)=  \mathbf{W}^1,$  а
  $ \mathcal{N}_d\mathbf{v}= \nabla\text{div}\textbf{v} =  \nabla\text{div}  \nabla h$  при	
 	  	  $\mathbf{v}=\nabla h \ \in  \mathbf{A}^{2}$.
    
Согласно  оценкам  (1.12) и (1.13) при $s=0$  пространство
    $ \mathbf{W}^1\subset \mathbf{H}^{1}$ \,и  $\mathbf{A}^{2}\subset \mathbf{H}^{2}$.
    Пространство   $\mathbf{C}_0^{\infty}(G) \cap \mathbf{V}^0$ плотно в $\mathbf{V}^0$ и содержится в $ \mathbf{W}^1$;   следовательно,  $ \mathbf{W}^1$  плотно в $\mathbf{V}^0$.    
Аналогично,   $\mathbf{C}_0^{\infty} (G) \cap  \mathbf{A}^0_{\gamma}$  плотно в $ \mathbf{A}^0_{\gamma} $ и содержится в $ \mathbf{A}^{2}$;
следовательно, $\mathbf{A}^{2}$ плотно в $ \mathbf{A}^{0}$.

 Если  поля $\mathbf{u}$ и $\mathbf{v}$  в (1.14) принадлежат  $\mathcal{D}(S)$, то $\gamma  (\mathbf{n}\cdot \mathbf{u})=\gamma  (\mathbf{n}\cdot \text{rot}\,\mathbf{u})= 0$, 
 $\gamma  (\mathbf{n}\cdot \mathbf{v})= \gamma  (\mathbf{n}\cdot \text{rot}\,\mathbf{v})=0$ ,  интеграл по $\Gamma$  зануляется \cite{yogi}  и  это равенсто принимает   вид:
$(\text{S} \mathbf{u}, \mathbf{v})= (\mathbf{u},\text{S}\mathbf{v})$. 

 Аналогично, если  поля  $\mathbf{u}=\nabla g$ и $\mathbf{v}=\nabla h$ в (1.15) принадлежат  $\mathcal{D}(\mathcal{N}_d)$ , то  $\gamma  (\mathbf{n}\cdot \mathbf{u})\equiv  \gamma\,(\mathbf{n}\cdot \,\nabla) g=0$,  $ \gamma  (\mathbf{n}\cdot \mathbf{v})\equiv   \gamma\,(\mathbf{n}\cdot \,\nabla) h=0$, интеграл по $\Gamma$  равен нулю и  это равенсто принимает   вид:
$
({\mathcal{N}_d} \mathbf{u}, \mathbf{v})= (\mathbf{u},{\mathcal{N}_d}\mathbf{v})$. 

Доказано, что $S$ и  $ \mathcal{N}_d$ -  самосопряженные расширения     операторов  $\text{rot}$ и  $\nabla\text{div}$ в  $\mathbf{L}_{2}(G)$
(см.  \cite{yogi,ds18}).
   \subsection{Гладкость собственных полей   операторов  $ \text{rot}$ и  $\nabla\text{div}$} 
         Спектральные задачи  для операторов  $ \text{rot}$ и  $\nabla\text{div}$ состоят в нахождении ненулевых полей $\mathbf{u}$ и $\mathbf{v}$ и чисел $\lambda$ и $\mu$  таких, что
\begin{equation*}
\klabel{srot__1_}
\text{rot} \mathbf{u}=\lambda \mathbf{u}(\mathbf{x}), \quad
\mathbf{x}\in G, \quad \gamma \mathbf{n}\cdot \mathbf{u}=0,
\quad \mathbf{u}\in \mathbf{C}^1(G)\cap \mathbf{C}(\overline{G}),
\eqno{(1.17)}
\end{equation*}
\begin{equation*}
\klabel{sgrd__1_}
\nabla\text{div}\mathbf{v}=\mu \mathbf{v}(\mathbf{x}), \quad
\mathbf{x}\in G, \quad \gamma \mathbf{n}\cdot \mathbf{v}=0,
\quad \mathbf{v}\in \mathbf{C}^2(G)\cap \mathbf{C}(\overline{G}).
\eqno{(1.18)}\end{equation*} 
Из Теорем 1,2  вытекают  важные свойства
решений  спектральных задач операторов  {\it  ротор и градиент дивергенции}: 

a)   {\it каждое ненулевое   собственное значение имеет  конечную кратность},

 b)  в  области $G$ с  гладкой границей их  обобщенные собственые поля   из $\mathbf{L}_{2}(G)
 $  является гладкими вплоть до границы. 

Доказательство.  Пусть  $\lambda\neq 0$,  а   поле $\mathbf{u}(\mathbf{x})$- соотвествующее ему решение задачи (1.17).  Такое поле $\mathbf{u}(\mathbf{x})$   есть решение однородной эллиптической задачи:
\[ \text{rot} \mathbf{u}=\lambda \mathbf{u}(\mathbf{x}), \quad \lambda\,\text{div} \mathbf{u}(\mathbf{x})=0, \quad
\mathbf{x}\in G, \quad \gamma \mathbf{n}\cdot \mathbf{u}=0,
\quad \mathbf{u}\in \mathbf{C}^1(G)\cap \mathbf{C}(\overline{G}).
\eqno{(1.19)}\]
Согласно Теореме 1 эта задача имеет конечное число линейно независимых решений   $\mathbf{u}_1(\mathbf{x}), ...,  \mathbf{u}_l(\mathbf{x})$, где  $l$ зависит от $\lambda$ и не зависит от  $\mathbf{u}$.  
Утверждение a) доказано. 

Решение   $\mathbf{u}(\mathbf{x})$  задачи (1.19)  принадлежит $ \mathbf{L}_{2}(G)$,  так как    $\|\mathbf{u}\|^2\equiv \int _G (\mathbf{u} \cdot \mathbf{u})d\mathbf{x} \leq V max_{\overline{G}}|\mathbf{u} \cdot \mathbf{u}|= V \|\mathbf{u} \cdot \mathbf{u}\|_{C(\overline{G})}$, где постоянная $V=\int_G 1\,d\mathbf{x}$.  

 Согласно (1.19)  $ \text{rot} \mathbf{u}_j=\lambda \mathbf{u}_j, \,\, \text{div} \mathbf{u}_j =0$ в $ G$, 
  $\gamma \mathbf{n}\cdot \mathbf{u}_j=0$.  Поэтому
$ \|\text{rot} \mathbf{u}_j \|= |\lambda| \| \mathbf{u}_j \|$ и
   оценка (1.12)  при   $s=0$ принимает вид: 
   $ C_0\|\mathbf{u}\|_{1} \leq (|\lambda|+1)  \|\mathbf{u}\|_{0} $, причём  постоянная $ C_0>0$. Значит,    $\mathbf{u}_j(\mathbf{x})$ принадлежит $\mathbf{H}^{1}(G)$ и
    \[ \|\mathbf{u}_j\|_{1} \leq  C_0 ^{-1}(|\lambda|+1)  \|\mathbf{u}_j
  \|_{0}, \quad    \|\mathbf{u}_j\|_0\leq \sqrt{V} \|\mathbf{u}_j\cdot \mathbf{u}_j  \|^{1/2}_{  {C}(\overline{G})}.   \eqno{(1.20)}\]
 Далее,  пусть   $s>0$ целое. Так как 
  $\| \text{rot}\mathbf{u}_j\|_s =|\lambda|\|\mathbf{u}_j(\mathbf{x})\|_s$,
 из оценки (1.12)   по индукции получаем:   
  $$\|\mathbf{u}_j\|_{s+1} \leq  C_s ^{-1} (|\lambda|+1)  \|\mathbf{u}_j\|_{s}\leq ...\leq  C_s ^{-1} ... C_0 ^{-1} (|\lambda|+1)^s  \|\mathbf{u}_j\|_{0}.             \eqno{(1.21)}$$

   Значит, поле   $\mathbf{u}_j(\mathbf{x})$ принадлежит $\mathbf{H}^{s+1}(G)$ для любого  целого  $s\geq 0$.
   
     {\it  Замечание 1.}  Известны вложения  пространств $H^{l+2}(\Omega)$  в $C^{l} (\bar {\Omega})$ при $l \geq  0$   в трехмерной области $\Omega$ и оценка 
    $ \|g\|_{C^{l} (\bar{\Omega})}\leq c_l \|g\|_{H^{l+2} (\Omega)}$ для любой функции   $g\in H^{l+2} (\Omega)$, причем постоянная 
    $c_l  > 0$ не зависит от $g$  [2, Теорема 3, § 6.2]. 
   
   Итак, поля  $\mathbf{u}_j(\mathbf{x})$ принадлежат $ \mathbf {C}^{l} (\bar G)$ для любого  целого  $l\geq 0$. Уверждение    b)  для ротора доказано.

    Аналогично, при  $\mu\neq 0$   собственное поле $\mathbf{v}(\mathbf{x})$  оператора $\nabla\text{div}$  есть решение однородной эллиптической задачи:
\[\nabla\text{div}\mathbf{v}=\mu \mathbf{v}(\mathbf{x}), \quad
\text{rot}\,\mathbf{v}=0, \quad\mathbf{x}\in G, \quad \gamma \mathbf{n}\cdot \mathbf{v}=0,
\quad \mathbf{v}\in \mathbf{C}^2(G)\cap \mathbf{C}(\overline{G}).
\eqno{(1.22)}\]
Согласно Теореме 2 эта задача имеет конечное число линейно независимых решений   $\mathbf{v}_1(\mathbf{x}), ...,  \mathbf{v}_k(\mathbf{x})$, где  $k$ зависит от $\mu$ и не зависит от  $\mathbf{v}$.  Утверждение a) доказано. \quad   
Любое решение   $\mathbf{v}_j(\mathbf{x})$  задачи (1.22)  принадлежит $\mathbf{L}_{2}(G)$,  так как   $\|\mathbf{v}\|^2_{\mathbf{L}_{2}(G)}\leq V\|\mathbf{v} \cdot \mathbf{v} \|_{  {C}(\overline{G})}$, где постоянная $V=\int_G 1\,d\mathbf{x}$.

 Ввиду того, что
   $ \|\nabla\text{div}\mathbf{v}\|= |\mu| \|\mathbf{v}\|$   в  $\mathbf{L}_{2}(G)$,
оценка (1.13)  при   $s=0$ принимает вид: 
$ C_0\|\mathbf{v}\|_{2} \leq (|\mu|+1)  \|\mathbf{v}\|_{0} $, причём  постоянная $ C_0>0$.

 Значит,    $\mathbf{v}_j(\mathbf{x})$ принадлежит $\mathbf{H}^{2}(G)$, и
\[ \|\mathbf{v}_j\|_{2} \leq  C_0 ^{-1}(|\mu|+1)  \|\mathbf{v}_j
\|_{0}, \quad    \|\mathbf{v}_j\|_0\leq \sqrt{V} \|\mathbf{v}_j\cdot \mathbf{v}_j  \|^{1/2}_{  {C}(\overline{G})}.    \eqno{(1.23)}\]
Далее,  пусть   $s>0$ целое. Так как 
$\|\nabla\text{div}\mathbf{v}_j\|_s =|\mu|\|\mathbf{v}_j(\mathbf{x})\|_s$,
из оценки (1.13)   по индукции получаем:   
$$\|\mathbf{v}\|_{2s+2} \leq  C_{2s} ^{-1} (|\mu|+1)  \|\mathbf{u}\|_{2s}\leq ...\leq  C_{2s} ^{-1} ... C_0 ^{-1} (|\mu|+1)^s  \|\mathbf{u}\|_{0}.      \eqno{(1.24)} $$

Значит,    $\mathbf{v}_j(\mathbf{x})$ принадлежит $\mathbf{H}^{2s+2}(G)\subset \mathbf{C}^{2s} (\bar {G}) $ для любого  целого  $s\geq 0$.
Утверждение b) доказано.

 \subsection{Гладкость базисных  полей   пространств $\mathcal{A}_H$ и $\mathcal{B}_H$}

Пространства   $\mathcal{A}_H$ и  $\mathcal{B}_H$  определяются решениями эллиптических систем (1.4)  и  (1.5)  в $ \mathbf{L}_{2}(G)$. Из формул (1.8) видно, что  компоненты  этих решений являются гармоническими функциями,  а значит,  они имеют непрерывные производные любого порядка.  Это впервые заметил  Герман Вейль для решений системы  (1.5) (\cite{hw} Теорема 1).
Краевые задачи   (1.4)  и  (1.5) удовлетворяют  условиям    В.Солонникова в Теореме 1.1  работы \cite{so71}. Откуда
получаем, что  пространства   $\mathcal{A}_H$ и  $\mathcal{B}_H$  конечномерны и их базисные поля  $\mathbf{g}_i(\mathbf{x})$ и $\mathbf{h}_j(\mathbf{x})\in 
\mathbf{C}^{\infty} (\bar G)$,  $i=1,..,\rho_1<\infty $, $j=1,..,\rho<\infty$.    
Отметим, что для $\mathbf{g}_i(\mathbf{x})$ и $\mathbf{h}_j(\mathbf{x})$ имеются оценки вида  (1.21) с   $\lambda=0 $.   
Borchers W.,  Sohr  H.  доказали \cite{boso}, что число  $\rho$ есть род границы $\Gamma$ области  $G$. 
В частности, если область  $\Omega$  гомеоморфна шару,  то $\rho=0$.


Если область  $\Omega$  гомеоморфна шару,  а $\mathbf{u}$  -решение задачи  $(1.5)$,    определяющей   $\mathcal{B}_H$, то   $\mathbf{u}= \nabla\, h$,  а функция $h$ - решение задачи Неймана для оператора Лапласа: $\Delta\,h =0 $ в $\Omega$,  $ \gamma (\mathbf{n}\cdot \nabla)\, h=0$.
Решение этой задачи   $N$ есть произвольная постоянная  $h=Const$. 
Следовательно,   $\mathbf{u}\equiv 0$ и пространство  $\mathcal{B}_H$ пусто.

Решение задачи (1.25)  в шаре $B, |\mathbf{x}|<R$,  сводится к задаче: $\Delta\,h =C$ в $B$,  где $C$ -произвольная постоянная,  с условием  Неймана      $\gamma (\mathbf{n}\cdot \nabla)\, h=0$.  Пусть $\mathbf{r}=\mathbf{x}$ - радиус-вектор, тогда нормаль  $\mathbf{n}=\mathbf{x}/r$   на границе шара, где
  $r= |\mathbf{x}|=R$.  Частное решение уравнения Пуассона  $\Delta\,h =C$ имеет вид: $h= 1/6C |\mathbf{x}|^2= 1/6Cr^2 $.
  Дифференцируя  по $r$,  получаем  $\gamma (\mathbf{r}\cdot \nabla)\, h= 1/3Cr|_{r=R}= CR/3$. Граничное условие Неймана принимает вид: $ CR/3=0$.  Значит,  $ C=0$ и  пространство 
 $\mathcal{A}_H(B)$ в шаре $B$ пусто.

 \subsection{Ортогональные базисы  в  $\mathcal{A}$, \, $\mathcal{B}$\,  и   в  $\mathbf{L}_{2}(G)$} Пространство $\mathbf{A}^2$ плотно в $\mathbf{A}^0 $ и
$\mathbf{A}^2\subset \mathbf{H}^2$.
Собственные поля $\mathbf{q}_{j}(\mathbf{x})$  оператора  $\nabla\mathrm{div}$ с ненулевыми собственными значениями ${\mu}_{j}$  принадлежат $\mathbf{A}^2$.

 {\small Множество собственных значений $\mu=-\nu^2$ этого оператора 	счётно, отрицательно и каждое из них имеет
	конечную кратность. Перенумеруем их в порядке возрастания их модуля: 	$0<-\mu_1\leq -\mu_2\leq ...$,  повторяя $ \mu_k$ столько раз, какова	его кратность. Соотвествующие вектор-функции  обозначим через
	$\mathbf{v}_{1}, \mathbf{v}_{2}$, ..., так чтобы каждому
	значению $\mu_{k}=-\nu^2_k$ соответствовала только
	одна  функция $\mathbf{v}_{k}$: $\nabla \mathrm{div}
	\mathbf{v}_{k}=-\nu^2_k \mathbf{v}_{k}$,  $ \gamma\mathbf{n}\cdot\mathbf{v}_{k}=0, $
	$k=1,2,...$.	
	Собственные функции, соответствующие одному и тому же
	собственному значению, выберем ортонормальными, используя процесс
	ортогонализации Шмидта  (см. \cite{vla}). 
	Поля, соответствующие
	различным с.- значениям, ортогональны. Их нормируем.
	Нормированные собственные поля  градента дивергенции обозначим 
	$\mathbf{q}_{l}$,   \quad  $ l=1,2,...$,  норма  $\|\mathbf{q}_{l}\|=1$.
	Они составляют полный ортонормированный базис в классе 
$\mathbf{A}^{0}$. 
  Зафиксируем его.}

Аналогично строится базис в классе $\mathbf{V}^{0}$ 
 \cite{saVS20}.
 
 {\it Замечание.}  Согласно (1.9)  оператор $\Delta \mathbf{u} \equiv -\mathrm{rot}^2 \mathbf{u}$  при  $\mathbf{u}\in \mathcal{B}$.   
 Собственные векторы ротора    всегда  встречаются парами: 
  каждому с.-полю    $\mathbf{u}^{+}_{j}$ с  $\lambda_j>0$  соответствует с.-поле   $\mathbf{u}^{-}_{j}$ с  $-\lambda_j$.  	
Это их свойство  в \cite{yogi} не отмечено.

Зафиксируем в  $\mathbf{V}^{0}$ ортонормированный базис $\{\mathbf{q}^{+}_{j}, \mathbf{q}_{j}^{-}\},  \quad  \|\mathbf{q}^{\pm}_{j}\|=1$: 
\[\mathrm{rot}
\mathbf{q}_{j}^{\pm}=\pm\lambda_j\, \mathbf{q}_{j}^{\pm}, \quad   \gamma\mathbf{n}\cdot\mathbf{q}_{j}^{\pm}=0,   \quad   j=1, 2, ...,  \quad  \mathbf{q}^{\pm}_{j}\in \mathbf{C}^\infty(\bar{G}).     \eqno{(1.26)} \]
Учитывая базисы пространств $\mathcal{A}_H$, \, $\mathcal{B}_H$
согласно (1.3), (1.6) видим, что  объединение $\{g_l\}$,	$\{\mathbf{q}_{l}\}$,  $\{h_j\}$  и $\{\mathbf{q}^{+}_{j}, \mathbf{q}_{j}^{-}\}$ есть базис объемлющего пространства  $\mathbf{L}_{2}(G)$.

     \subsection{Явный вид собственных полей   ротора в шаре  $B$}  	
    Спектральные задачи для операторов      	ротор и градиент дивергенции  в шаре    решены автором   полностью  в   \cite{saUMJ13}. 
   
    Имеется несколько способов   решения  первой задачи  \cite{wo58, cdtgt,saUMJ13}.
   
   Учитывая  приложения \cite{wo58} и конкурирующие интересы           \cite{ cdtgt},  кратко изложим наш путь  решения   этой задачи  \cite{saUMJ13}.

  	Собственные числа $\lambda_{n,m}$  ротора  в шаре радиуса $R$ равны $\pm \rho_{n,m}/R$,  где числа $\pm \rho_{n,m}$ - нули функций  $\psi_n(r)$:
  	\begin{equation*}\label{bes  2}
  	\psi_n(z)
  	= (-z)^n\left(\frac d{zdz}\right)^n\left(\frac{\sin z}z\right),  \quad  m,n\in {\mathbb {N}}. \eqno{(1.27)}
  	\end{equation*}
  Функции  $\psi_n(z)$- это цилиндрические  функции $J_{n+1/2}(z)$,  где  $n\geq 0$ целое.  Это заметил ещё Леорнард Эйлер (см. \cite{vla} $\S23$).
  
  Числа $\pm \rho_{n,m}$ и  $ \rho^2_{n,m}>0$ - нули функций  $\psi_n(z)$.
  
  Кратность  собственного значения  $\lambda^{\pm}_{n,m}$ равна $2n+1$. 
  
  Пусть $\mathbf{i}_r, \mathbf{i}_\theta,
  \mathbf{i}_\varphi$-репер, поле  $\mathbf{u}=u_r\,\mathbf{i}_r+ u_{\theta} \mathbf{i}_{\theta}+ u_{\varphi}
  \mathbf{i}_\varphi$.

  	{\bf Формулы решений задачи (1.26).}
    		 {\it  Ненормированные собственные поля
    		 	 $ \mathbf{u}_{\kappa
  		}^{\pm }$ задачи  (1.26) в сферических координатах вычисляются по
  		формулам:}
  		\begin{equation*}
  		\label{vrsff   2}\begin{array}{c}
   \mathbf{u}_{\kappa}^{\pm }=c_{\kappa }^{\pm }(\pm\lambda _{n,m}
  		r)^{-1}{\psi }_{n}
  		(\pm\lambda _{n,m}r)Y_{n}^{k}(\theta ,\varphi  )\,\mathbf{i}_r+\\
  		c_{\kappa }^{\pm }{{(\pm\lambda _{n,m}^{\pm }r)}^{-1}}
  		Re[\Phi _{n}(\pm\lambda _{n,m}r)](Re HY_{n}^{k}\,
  		\mathbf{i}_\varphi+
  		Im HY_{n}^{k}\,\mathbf{i}_\theta)+\\
  		c_{\kappa }^{\pm }{{(\pm\lambda _{n,m}r)}^{-1}}
  		Im[\Phi _{n}(\pm\lambda _{n,m}r)](-Im HY_{n}^{k}\,
  		\mathbf{i}_\varphi+
  		Re HY_{n}^{k}\,\mathbf{i}_\theta).\end{array}   \eqno{(1.28)}\end{equation*}
  		где
  		$Y_{n}^{k}(\theta ,\varphi  )$--сферические функции, числа
  	$c_{\kappa }^{\pm }\in \mathbb{R}$ -произвольны, $\kappa=(n,m,k)$- мульти-индекс,
  $m{{,}^{{}}}n\in \mathbb{N}$,    $|k|\le n$, 
   	\[\Phi _{n}(\lambda\,r)=
  	\overset{\mathop{{}}}\,\int\limits_{0}^{r}{}\,{{e}^{i\lambda   			(r-t)}}^{{}}
  	{\psi }_{n}(\lambda t){t}^{-1}dt,
  	\quad Im \Phi _{n}(\pm\rho_{n,m})=0,        \eqno{(1.29)}\]
  	\begin{equation*}
  	\label{oph   1}
  	\text{H}v=
  	\left(\sin ^{-1}\theta {\partial }_{\varphi }+
  			i{\partial }_{\theta } \right)v,  \quad 	\text{K}w=
  		\sin ^{-1}\theta 	\left({\partial }_{\theta } \sin \theta+
  			i{\partial }_{\varphi } \right)w .   \eqno{(1.30)} \end{equation*}  
  Решению  	этой	   спектральной задачи
 способствовали  наблюдения автора:

    1. {\it 	Пусть поле $\mathbf{u}$ -  решение  спектральной задачи (1.19) в шаре $B$,  	$\mathbf{x}$-радиус- вектор, а
  	  $v(\mathbf{x})$- их  скалярное  		произведение   $\mathbf{x}\cdot \mathbf{u}=r\,u_r$.
        Тогда	функция    $v(\mathbf{x})$  есть решение  спектральной задачи Дирихле для уравнения Лапласа: }
  	\[	\label{ldo__1_}           -\Delta v=\lambda^{2}\,v  \quad \text{в} \quad B,
  \quad v|_{S}=0,   \quad \text{с условием} \quad  v(0)=0.                 \eqno{(1.31)}\]

  		2. {\it  Уравнения $\mathrm{rot}\mathbf{u}=\lambda \mathbf{u}, \,  \mathrm{div}\mathbf{u}=0$, записанные   в сферических координатах,  представляются 	в виде двух комплексных уравнений   			
  			\[({\partial}_r - i\lambda) r\,w= r^{-1}H\,v,
  			\quad  K\,w=\lambda\,v-i\,r^{-1\,}\,{\partial}_r (r\,v),    \eqno{(1.32)} \]
  			относительно функций $v=ru_r$ и $w=u_{\varphi}+iu_{\theta}$	 с операторами  $H$  и $K$  в (1.30).  
  		 }
  	 
  	 3.  {\it   Уравнения (1.31) на  функцию $v$ являются условиями   	 совместности системы (1.32).}
  
  		Таким образом,  решение задачи сводится к решению:
  		  		
  	{ \small $1^0)$  спектральной задачи Дирихле - Лапласа  (1.31).  Её решения- пары  $ \lambda_{\kappa}^2=(\rho_{n,m}/R)^2$
  		и  ${v}_{\kappa}=c_{\kappa }{\psi }_{n}(\rho_{n,m}r/R )Y_{n}^{k}(\theta ,\varphi  )$, такие что ${\psi }_{n} (\rho^2 _{n,m})=0$
  		. Условие 	  v(0)=0    обеспечивается обнулением постоянных $c_{\kappa }=0$  при $\kappa=(0,m,0)$ (см. В.С.Владимилов   \cite{vla} гл.V $\S26$).  Они определяют  $ \lambda^{\pm}_{\kappa}=\pm \rho_{n,m}/R$-собственные  значения задачи (1.17) в шаре $B$  и  функции    $ u_{r,\kappa }=v_{\kappa}/r $-  радиальные компоненты собственных полей. }
  	
  	{ \small $2^0)$ к интегрироваию уравнений (1.32) с  $ \lambda= \lambda^+_{\kappa}>0$ и $v=v^+_{\kappa}$, а затем с  $ \lambda= \lambda^-_{\kappa}<0$ и $v=v^-_{\kappa}$,  и  вычилению комплексных функций $w^{\pm}_{\kappa}$,     задающих  касательные компоненты полей   $\mathbf{u}^{\pm}_{\kappa}$; они 
  		определятся однозначно  условием: $w^{\pm}_{\kappa}\in L_2(B)$
  		
  		\small $3^0)$ к 	построению полей 	$\mathbf{u}_{\kappa}^{\pm}(\mathbf{x})\in \mathbf{L}_2(B)$. }    	
  	
  	В итоге, получаем список решений (1.28).
  	
 {\it Замечание.}  	Позже 	уравнение (1.31) на функцию  $v=r\,u_r$  при  минимальном  собственном значении $\lambda=4.4934.../R$ автор обнаружил  в статье \cite{cdtgt}.

   \subsection{Явный вид собственных полей    $\nabla\mathrm{div}$ в шаре  $B$}      Собственные значения  оператора $\nabla\mathrm{div}$ равны $-\nu_{n,m}^2$, где $\nu_{n,m}=\alpha_{n,m}/R$,  а  числа  $\alpha_{n,m}$ - нули производных  $\psi'_n(r)$,   $n \geq 0,\, m\in {\mathbb {N}}$ ;     	
 	кратность  собственных значений  $-\nu^2_{n,m}$    равна $2n+1$.	 	
 
 	Собственные поля $\mathbf{v}_{\kappa}= \nabla  g_{\kappa}$    градиента дивергенции  - решения задачи:
 	\[   \nabla \mathrm{div}
 	\mathbf{v}_{k}=-\nu^2_{\kappa}\mathbf{v}_{\kappa}, \quad  \gamma\mathbf{n}\cdot\mathbf{v}_{\kappa}=0,  \quad  \mathbf{v}_{\kappa}= \nabla  g_{\kappa}\in \mathcal{C}^\infty(\bar{G}).   \eqno{(1.33)} \] 
Так как   $ \nabla \mathrm{div}\nabla  g_{\kappa}\equiv  \nabla 
\Delta_c\,g_{\kappa}= 
\Delta_c \, (\nabla g_{\kappa}) =  -\nu^2_k \,(\nabla g_{\kappa}),  \quad \gamma(\mathbf{n}\cdot  \nabla ) g_{\kappa}=0$ эта задача сводится к задаче Неймана для скалярного оператора Лапласа и градиенту фукций $g_{\kappa}$.  Матричный $(3 \times 1)$ оператор   $ \nabla \mathrm{div}\nabla  g \equiv  \nabla 
\Delta_c\,g$ эллиптичен.

Соответствующие  $-\nu^2_{\kappa}\equiv  -\nu^2_{n,m}$ собственные функции  $g_{\kappa}$ имеют вид:
\[ g_{\kappa }(r, \theta ,\varphi)=c_{\kappa }\psi _{n}
(\alpha _{n,m}r/R)Y_{n}^{k}(\theta ,\varphi  )    \eqno{(1.34)} \]

Поля   $\mathbf{v}_{\kappa}= \nabla  g_{\kappa}$  являются решениями задачи (1.25); их компоненты $(v_r, v_{\theta}, v_{\varphi })$  имеют вид:
 $ v_{r, \kappa }(r, \theta ,\varphi)=c_{\kappa } (\alpha _{n,m}/R)\psi'_{n} (\alpha _{n,m}r/R)Y_{n}^{k}(\theta ,\varphi), $ 
 \[ (v_{\varphi}+iv_{\theta})_{ \kappa}=c_{\kappa } (1/r)       \psi_{n} (\alpha _{n,m}r/R)\,H\,Y_{n}^{k}(\theta ,\varphi  )  \eqno{(1.35)}                  \]
 
 При $ \kappa=(0,m,0)$ функция $Y_{0}^{0}(\theta ,\varphi)=1, HY_{0}^{0}(\theta ,\varphi)=0,$  поэтому 
 
\[ v_{r, (0,m,0) }(r)=  c_{ (0,m,0)}(\alpha _{0, m}/R) \psi'_{0} (\alpha _{0,m}r/R) , \quad  (v_{\varphi}+iv_{\theta})_{ (0, m, 0)}=0.       \eqno{(1.36)}             \]

  Построенный в  шаре  $B$  базис из собственных полей операторов градиента дивергенции и  ротора является полным в $\mathbf{L}_{2}(B)$,
		так как $\mathbf{L}_{2}(B)=\mathbf{A}^0 \oplus \mathbf{V}^0$.

 \subsection{Потоки с  минимальпой энергией,  визуализация  их}
 	Эти формулы используются при рассчетах поля скоростей $\mathbf {u}_{\kappa }^{\pm}(\mathbf {x})$ и визуализации вихревых потоков.  
 Формулы  полей  $\mathbf{u}_{\kappa}^{\pm}(\mathbf{x})$  при  $n=1$, $\kappa=(1,1,0)$ и  $\kappa=(1,1,\pm 1)$ выражаются  наиболее просто. Так,  компоненты поля $\mathbf{u}^{+}_{(1,1,0)}(\mathbf{x})$ имеют вид:    
 \[u_r=2\rho( r\rho)^{-3}( sin\,(r\rho)- r\rho cos\,(r\rho)) cos\,\theta,   \] 
 \[u_\theta =( r\rho)^{-3}( sin\,(r\rho)- r\rho \,cos\,(r\rho)-( r\rho)^{2}sin\,(r\rho)) sin\, \theta ,       \eqno{(1.37)} \]
 \[u_{\varphi} = ( r\rho)^{-2}(( sin\,(r\rho)- r\rho \,cos\,(r\rho)) sin\, \theta.  \]  
 	Профессор Исламов Г.Г.\cite{gais}, используя эти формулы   и  программу  Wolfram Mathematica  осуществил
 визуализацию линий тока поля   $u_{1,1,0}^{+}(\mathbf {x})$ ротора  радиуса 1 со значением $\rho =\rho_{1,1}$ \footnote{ 
 	http://www.wolfram.com/events/ technology-conf.-ru/  2016/resources.html }.   
 Траектория движения  трёх соседних точек напоминает ленту, которая 	обматывает тороидальную катушку  (см.  катушка Исламова  в \cite{saVS20}).
  \footnote{Исламов Галимзян Газизович (02.02.1948-22.11.2017) }
  
 	В связи с задачами  астрофизики  S. Chandrasekhar, P.С. Kendall изучали  собственные поля ротора в шаре \cite{chake}  и  в цилиндре. 
 	 Они нашли   элементарный способ их вычисления  в цилиндре  (с условием периодичности вдоль оси).
 	
 	D.Montgomery, L.Turner, G.Vahala,  	 изучая магнито гидродинамическую турбулентность в цилиндре  \cite{motva},
 	 исползовали эти формулы.
 	 
 	  Они показали, что три интегральных инварианта имеют простые квадратичные выражения в терминах спектральных разложений.

 	J. Cantarella,  D. De Turck, H. Gluck and  M.Teitel  (лаб. 	 "Физика плазмы")   исследовали собственные поля ротора  в шаре  радиуса $b$  и в шаровом слое.
 
 Уравнение (1.31) на функцию  $v=r\,u_r$  при минимальном   собственном значении $\lambda_{1,1}=\rho_{1,1}/b>0$ автор обнаружил  в их статье \cite{cdtgt}.
 Они приводят  также   соответствующую	 $\lambda_{1,1}$  формулу  собственного поля ротора в шаре  (см. Theorem A). К сожалению,  с опечаткой: $1/\lambda$  вместо $\lambda$.  
 
 Исправив её,   мы получили   формулы  (1.37)  компонент поля  $\mathbf{u}^{+}_{(1,1,0)}(\mathbf{x})$. 

В Fig.1  в  \cite{cdtgt}  представлены 
интегральные кривые  поля $\mathbf{u}^{+}_{(1,1,0)}(\mathbf{x})$    и даётся их описание. \quad  Цитирую: "они заполняют семейство концентрированных "торов"  с  замкнутой орбитой  "ядра,"  \, типичных  для осесимметричных собственных полей ротора; специальная орбита  начинается на южном полюсе сферы в момент  времени $-\infty $,  проходит вертикально вверх по оси z и достигает северного полюса ко времени  $+\infty$; орбиты на граничной сфере   начинаются на северном полюсе в момент времени  $-\infty $,  продолжаются по линиям долготы к южному полюсу до момента времени  $+\infty $; имеются две стационарные точки в её полюсах."

В  статье \cite{cdtgt} отмечено, что Woltjer использовал  векторное 
поле   $\mathbf{u}^{+}_{(1,1,0)}(\mathbf{x})$  для моделирования магнитного поля в Крабовидной туманности  \cite{wo58}.

Мы провели независимое исследование этого поля. Галимзян Исламов демонстрировал эти орбиты  вживую на  экране в  МГУ  во время нашего совместного доклада на конференции  Бицадзе-100,  факультет ВМК 2016.

 Авторы  \cite{cdtgt}  приводят также формулы   базисных  полей ротора   в шаре  для других собственных значений.  
 Для сравнения  мы  приводим  весь список  (1.28).

   \subsection{Степени оператора Лапласа в классах $\mathcal{A}$, \, $\mathcal{B}$\,  и   в  $\mathbf{L}_{2}(G)$}
Из формул (1.9) при $k=2, 3, ...$ имеем
\[\Delta^k \mathbf{v}\equiv  (\nabla \mathrm{div})^k \mathbf{v}   \quad  \text{при}    \quad  \mathbf{v}\in \mathcal{A},  \quad  \Delta^k \mathbf{u} \equiv (-1)^k(\mathrm{rot})^{2k}\,\mathbf{u}  \quad  \text{при}  \quad      \mathbf{u}\in \mathcal{B}.   \eqno(1.38)\]
В  $\mathbf{L}_{2}(G)$ оператор $\Delta^k$
выражается через  $(\nabla \mathrm{div})^k$ и $(\mathrm{rot})^{2k}$,   а также через скалярный оператор $\Delta_c^k = (\partial_1^2+\partial_2^2+\partial_3^2)^k$ : 
\begin{equation*}\klabel{dd 1}\mathrm{\Delta}^k\,\mathbf {v}
= (\nabla \mathrm{div})^k\,\mathbf {v} + (-1)^k\,(\mathrm{rot})^{2k}\, \mathbf {v}= \Delta^k_c\, I_3\,\mathbf {v} , \quad \text{где}
\,\,  \mathbf{v}=(v_1, v_2, v_3). \eqno(1.39) \end{equation*} \quad  
Эти формулы следуют из формул (1.8),  учитывая, что операторы $\mathrm{rot}$ и $\nabla \mathrm{div}$ аннулируют друг друга.  Они являются проекторами :
 $\nabla\mathrm div$ проектирует $\mathbf{L}_{2}(G)$ на $\mathcal{A}$, а  $\mathrm{rot}$ -   на $\mathcal{B}$   .

С.Л.Соболев  изучил  периодическую задачу $\pi$ и краевые задачи $D$ и $N$ для скалярного полигармонического уравнения $\Delta^m\,u=\rho$ в пространствах ${W} _2^m(\Omega)$ c правой частью -- обобщённой функцией  (см. \cite{sob},   $\S 9$ гл. 12).

В периодическом случае, например,  он доказал    теорему  (цитирую):

Т е о р е м а  XII.13.  {\it Оператор $\Delta^m$ переводит произвольную функцию $u$  из   $\bar {W} _2^{(m)}$ в 
$\Delta^m\,u=\rho$ --  элемент  $\bar {L}_2^{(m)^*}$.
Обратно, для произвольной обобщённой функции $\rho$  из
 $\bar {L}_2^{(m)^*}$ существует   функция $u\in \bar {W} _2^{(m)}$  такая, что $\Delta^m\,u=\rho$. 
   Эта функция определяется  с точностью  до  произвольного постоянного слагаемого. }

Операторы   $(\nabla \mathrm {div})^p$ и $(\mathrm {rot})^{2q}$, где  $ p$ и  $ q$ - натуральные  числа,
- аналоги  полигармонических  операторов $\Delta^m$
 в классах $\mathcal{A}$   и  $\mathcal{B}$ (см. (1.38)).
Мы покажем,  что 
 оператор $(\nabla \mathrm {div})^{2p}$ переводит
произвольное поле $\mathbf{w}$  из   $ {A}^{2p}$ в 
$(\nabla \mathrm {div})^{2p}\,\mathbf{w}=\rho$ --  элемент  ${A}^{-2p}\equiv  ({A}_0^{2p})^*$;
а оператор $(\mathrm {rot})^{2q}$ переводит
 произвольное поле $\mathbf{u}$  из   $ {W}^{q}$ в 
$(\mathrm {rot})^{2q}\,\mathbf{u}=\mathbf{v}$ --  элемент  ${W}^{-q}\equiv  ({W}_0^{q})^*$.

Доказаны и обратные утверждения.

  \section{Пространство $\mathcal{A}$ потенциальных полей }

 В статье автора \cite{ds18} детально рассмотрена структура класса $\mathcal{A}$ потенцциальных полей, его базис    и оператор $ \mathcal{N}_d$.  Здесь мы рассмотрим  его подпространства  $\mathbf{A}^{2k}$.
    По  определению ${\mathcal{{A}}}(G) =\{\nabla h, h\in H^1\}$,  $\mathcal{A}_H $ - ядро оператора  $ \nabla\mathrm{div}$  в $\mathcal{A}$,  а    $\mathbf{A}^0$ - его ортогональное дополнение в  $\mathcal{A}$,  $\mathcal{A}=\mathcal{A}_H\oplus \mathbf{A}^0$,    \quad
     $ \mathcal{A}_{\gamma} (G) =\{\nabla h, h\in H^2(G): \gamma(\mathbf{n}\cdot \nabla) h=0\}$,   \quad
    $ \mathcal{A}^0_{\gamma}  = \mathbf{A}^0\cap \mathcal{A}_{\gamma}$
   
  Подпространство
$\mathbf{A}^2=\{ \mathbf{v}\in
\mathcal{A}^0_{\gamma}: \nabla\mathrm{div} \mathbf{v}\in
\mathcal{A}^0_{\gamma}\}$  есть область определегия оператора   $ \mathcal{N}_d$;  оно плотно в $\mathbf{A}^0$ и
$\mathbf{A}^2\subset \mathbf{H}^2$ (согласно п. 1.5). 
Собственные поля $\mathbf{q}_{j}(\mathbf{x})$  оператора  $\nabla\mathrm{div}$ с ненулевыми собственными значениями $(-\nu^2_j)$:
 $\nabla \text{div}\, \mathbf{q}_{j}= -\nu^2_j\,\mathbf{q}_{j}$, \,
$\gamma(\mathbf{n}\cdot \mathbf{q}_{j})=0$,
принадлежат пространству  
$\mathbf{A}^2$.
Они составляют ортонормальный базис $\{\mathbf{q}_{j}(\mathbf{x})\}$  в 	$\mathbf{A}^{0}$.  
  	Проекция поля $\mathbf{f}\in \mathbf{L}_2(G)$ на
$\mathbf{A}^0$ имеет вид:  \[\mathcal{P}_{\mathcal{A}}\mathbf{f}\equiv 	\mathbf{f}_{\mathcal{A}}(\mathbf{x})=	
\lim_{n\rightarrow\infty}(\mathbf{f}^n_{\mathcal{A}})=\sum_{j=1}^\infty
(\mathbf{f},\mathbf{q}_{j})\mathbf{q}_{j}(\mathbf{x}),   \eqno{(2.1)}\]	
где	 $\mathbf{f}^n_{\mathcal{A}}$--  частичные суммы  этого
ряда. 

Оператор $ \mathcal{N}_d$  определен 
  и совпадает с $\nabla\mathrm{div}$ на  $ \mathbf{A}^{2}$, поэтому	
\[ \mathcal{N}_d\mathbf{f}_{\mathcal{A}}=
\lim_{n\rightarrow\infty}
\nabla\mathrm{div}\,(\mathbf{f}^n_{\mathcal{A}})
=-\sum_{j=1}^{\infty}\nu^2_j
(\mathbf{f},\mathbf{q}^{}_{j})\mathbf{q}^{}_{j}(\mathbf{x}),   \eqno{(2.2)}\]
если ряд	сходится и принадлежит $\mathbf{A}^0$.	
Это так, если $f\in\mathbf{H}^2(G)$.  

Доказано, что 		оператор $\mathcal{N}_d$ замкнут и самосопряжён  \cite{ds18}.  

 \subsection{Подпространства  $\mathbf{A}^{2k}$ в $\mathcal{A}$}
Рассмотрим ещё пространства\footnote{Они совпадают с пространствами $\mathbf{A}^{2k}_ {\gamma}$  в  \cite{ds18}, если пространство  $\mathcal{A}_{H}$ пусто. }
\[\mathbf{A}^{2k}=\{\textbf{ f}\in  \mathcal{A}^0_ {\gamma},..., (\nabla  \text{div})^k\textbf{ f}\in \mathcal{A}^0_ {\gamma} \}, \quad k=1,2, ... \eqno{(2.3)}\]	  		  	  		  	 
{\it Замечание. }  Согласно  оценке (1.13) пространство $\mathbf{A}^{2k}\subset \mathbf{H}^{2k}$. 
Оно    является   проекцией пространства Соболева  $\mathbf{H}^{2k}$ порядка ${2k} $ на класс  $\mathcal{A}$, так как  для любого поля   $\mathbf{f}\in  \mathbf{H}^{2k} $  его проекция $\mathcal{P}_{A}\mathbf{f}\in  \mathbf{A}^{2k}$;
 если  же $\mathbf{f}\in   \mathbf{A}^{2k}$,   то $\mathcal{P}_{A}\mathbf{f}= \mathbf{f}$,  а 	его проекция 
  на 	  	  $\mathcal{B}$ равна 0.  

	  Пространство $\mathcal{A}^0_{\gamma}$ ортогонально ядру оператора
$\mathcal{N}_d$
в $\mathbf{L}_{2}(G)$, поэтому   $\mathcal{N}_d$ имеет единственный обратный	  	  оператор: 
\begin{equation*}\label{Nobr__2_}\mathcal{N}_d^{-1}\mathbf{f}_{\mathbf{A}}=
-\sum_{j=1}^{\infty}\nu_j^{-2}
(\mathbf{f},\mathbf{q}_{j})\mathbf{q}_{j}(\mathbf{x}).    \eqno{(2.4)}\end{equation*} 
	Оператор $\mathcal{N}_d^{-1}$ - компактен. 

Следствие. {\it Спектр оператора $\mathcal{N}_d^{-1}$ точечный с
	единственной точкой  накопления в нуле,\quad
	$\nu^{-2}_j\rightarrow 0$ при ${j\rightarrow\infty}$.}

\subsection{Сопряжённые пространства  $\mathbf{A}^{-2k}$}
По определению пространство 
 ${H}^{s}_0(G)$ есть замыкание	в норме ${H}^{s}(G)$  функций из ${C}^{\infty}_0(G)$.  
    $\mathcal{A}_0 =\{\nabla h, h\in H^1_0\}$, \quad
$\mathbf{A}^{2k}_ {0}=\{\textbf{ f}\in  \mathcal{A}_ {0}, ..., (\nabla   \text{div})^k\textbf{ f}\in \mathcal{A}_ {0} \}$.  \quad
 Пространство линейных непрерывных функционалов над   $\mathbf{A}^{2k}_ {0}$,  обозначим $(\mathbf{A}^{2k}_0)^*$. Они равны нулю на $ \mathcal{A}_H$ (см. п.2.3).
 
 В п. 2.4 мы покажем, что эти пространства можно отождествить с пространствами $\mathbf{A}^{-2k}$ 
  порядка $-2k$.
   Наконец, $ \mathcal{A}^*$- это объединение  $\mathbf{A}^{-2k}$ при  $k\geq 1$.
 
Цепь вложений пространств $\mathbf{A}^{2k}$ 
имеет вид: \[\subset \mathbf{A}^{2k}\subset...\subset\mathbf{A}^{2}\subset  \mathbf{A}^0\subset \mathbf{A}^{-2}\subset...\subset \mathbf{A}^{-2k}\subset \eqno{(2.5)}\]

Операторы $\mathcal{N}_d:\mathbf{A}^{2k} \rightarrow \mathbf{A}^{2(k-1)}$ \, обратимы  при    $k> 1$   и 
$$\|\mathcal{N}_d^{-1}\mathbf{f}\|^2_{\mathbf{A}^{2k}} 
\leq  c^2_k \|\mathbf{f}\|^2_{\mathbf{A}^{{2(k-1)}}},  \quad
\|\mathcal{N}_d  \mathbf{f}\|^2_{\mathbf{A}^{2(k-1)}}\leq  
c^{-2}_k \|\mathbf{f}\|^2_{\mathbf{A}^{2k}},       \eqno{(2.6)}$$
  где
$c^2_k= max _j(1+1/{\nu}_{j}^{2k})$,   а $1/{\nu}_{j}\to 0$ при $j\to \infty $. 

{\it Замечание.} Автор изучал также оператор  $ \mathcal{N}_d+\lambda I$ в  \cite{saVS20, ds18}, доказана

\begin{theorem} \klabel{Nd_3}Оператор $\mathcal{N}_d+\lambda I:\mathbf{A}^{{2(k+1)}} \rightarrow \mathbf{A}^{2k}$ -фредгольмов при $k\geq0$. 
	Если $\lambda \overline{\in} Sp (\mathcal{N}_d)$,  ,
	то 	 оператор  $ \mathcal{N}_d+\lambda I$  (и его обратный)  отображает    пространство  $ \mathbf{A}^{2(k+1)}$
	на  $  \mathbf{A}^{2k}$ (и обратно) 	взаимно	однозначно и непрерывно.\end{theorem}

Оператор     $ \mathcal{N}_d \mathbf{u}$ совпадает с $\nabla \mathrm{div}\mathbf{u}$, если   $\mathbf{u}\in \mathbf{A}^{2}\equiv \mathcal{D}(\mathcal{N}_d) $. Поэтому оператор $(\nabla \mathrm{div})^{k}$ на $\mathbf{A}^{2k}\subset \mathbf{A}^{2}$ совпадает с $ \mathcal{N}_d^{k}$ при    $k> 1$. 
 
\subsection{Оператор    $\mathcal{N}^{2k}_d$ в пространстве $\mathbf{A}^{2k}$}Основное утверждение.

  {\it Оператор    $\mathcal{N}^{2k}_d$ отображает пространство $\mathbf{A}^{2k}$  на  $\mathbf{A}^{-2k}$ и обратно.}

Этапы доказательства:

 Шаг 1-й:  {\it Оператор $ \mathcal{N}_d^{2k}$ отображает пространство
$ \mathbf{A}^{2k}$  на   $ (\mathbf{A}^{2k}_0)^*$.}

Действительно,  	пусть $\mathbf{w}$ произвольный элемент из $\mathbf{A}^{2k}$, а
$\mathbf{w}_{\eta}$ -- средняя вектор-функция для него,
$\mathbf{w}_{\eta}\in\mathbf{A}^{2k}_0$;   поле\  $\mathbf{u}\in \mathbf{A}^{2k}$.
Рассмотрим главную часть скалярного произведения в $\mathbf{A}^{2k}(G)$:	
$$(\mathbf{u},\mathbf{w}_{\eta})_{2k}\equiv 
((\nabla \mathrm{div})^{k}\,\mathbf{u},(\nabla \mathrm{div})^{k}\,\mathbf{w}_{\eta}).$$
Проинтегрируем по частям:	
\[	(\mathbf{u},\mathbf{w}_{\eta})_{2k}=
((\nabla \mathrm{div})^{2k}\,\mathbf{u},\,\mathbf{w}_{\eta})=
\int_G
\mathbf{v}\cdot (\mathbf{w}_{\eta})\, d\mathbf{x}.   \eqno{(2.7)}\]
Левая часть имеет
предел при $\eta\rightarrow 0$, равный $(\mathbf{u},\mathbf{w})_{2k}$.
Следовательно,  правая часть также будет иметь предел и интеграл
$\int_G \mathbf{v}\cdot \mathbf{w}\, d\mathbf{x}$ существует при
любой $\mathbf{w}\in\mathbf{A}^{2k}(G)$. Кроме того,  из     
неравенства Коши-Буняковского следует оценка этого интеграла:
$$\left|\int_G
\mathbf{v}\cdot \mathbf{w}\, d\mathbf{x}\right|\leq
\|\mathbf{u}\|_{\mathbf{A}^{2k}}\|\mathbf{w}\|_{\mathbf{A}^{2k}}.   $$ 
Значит, $\mathbf{v}$ есть линейный функционал из
$(\mathbf{A}_{0}^{2k})^*$.
 
Применим его к полям $\mathbf{g}_i$,
составляющим базис пространства $\mathcal{A}_H(G)$.  Учитывая,
что  $\nabla \mathrm{div}\,\mathbf{g}_i=0$,
получим \[\int_G \mathbf{v}\cdot \mathbf{g}_i\, d\mathbf{x}=0. \quad
i=1,...,\rho_1.                  \eqno{(2.8)}  \] Итак,  на  полях
$\mathbf{w}$, отличающихся на вектор-функцию $\mathbf{g}$ из
$\mathcal{A}_H(G)$, его значения совпадают.
Пусть 	$\mathcal{A}/\mathcal{A}_H$ - фактор-пространство   	$\mathcal{A}(G)$ по $\mathcal{A}_H$   (пространство классов смежности). $\mathcal{A}^{2k}(G)=\mathbf{A}^{2k}(G)/\mathcal{A}_H$, его элементы имеют вид: $\mathbf{w}+\mathbf{g}$, где   $\nabla \mathrm{div}\,\mathbf{g}=0$.

\subsection{Оператор    $\mathcal{N}^{2k}_d$ в фактор-пространстве	 $\mathcal{A}^{2k}$ }
	\qquad    Шаг 2-й:
	
{\it	 Пространство
$\mathcal{A}^{2k}(G)$ становится гильбертовым, если  ввести скалярное произведение}
$$\{\mathbf{u},\mathbf{w}\}_{2k}\equiv (\mathbf{u},\mathbf{w})_{2k}=
((\nabla \mathrm{div})^{k}\,\mathbf{u},(\nabla \mathrm{div})^{k}\,\mathbf{w}).   \eqno{(2.9)}$$

Воспользуемся ортонормированным базисом в $\mathbf{A}^{0}$.
При $ \mathbf{f}\in \mathcal{A}^{2k}$, 
$\mathbf{g}_\eta\in	\mathcal{A}^{2k}_0$, 	в	терминах рядов Фурье  оно имеет вид:
\[\{\mathbf{f},\mathbf{g}_\eta\}_{2k}\equiv 
( \mathcal{N}_d^{k}\,\mathbf{f}, \mathcal{N}_d^{k}\,\mathbf{g}_\eta)=
\sum_{j=1}^{\infty}{\nu}_{j}^{4k}
[(\mathbf{f},\mathbf{q}_{j})(\mathbf{g}_\eta,\mathbf{q}_{j})
],  \eqno{(2.10)}\]  так как  
\[ \mathcal{N}_d^k\mathbf{f}=
\lim_{n\rightarrow\infty}
(\nabla\mathrm{div})^k\,(\mathbf{f}^n_{\mathcal{A}})= (-1)^k
\sum_{j=1}^{\infty}\nu^{2k}_j
(\mathbf{f},\mathbf{q}^{}_{j})\mathbf{q}_{j}(\mathbf{x}),   \eqno{(2.11)}\]

Для того чтобы функционал $\rho$ служил элементом
$( \mathcal{A}^{2k}_0)^*$, нужно, чтобы скалярное произведение
$({\rho}(\mathbf{x}),\mathbf{w}(\mathbf{x}))$ существовало при всех
$\mathbf{w}(\mathbf{x})\in  \mathcal{A}^{2k}$ и удовлетворяло
неравенству: $({\rho}(\mathbf{x}),\mathbf{w}(\mathbf{x}))\leq
M_{2k}\,\| \mathcal{N}_d^k\,\mathbf{w}\|$. 	
Мы имеем
$$({\rho},\mathbf{w})=\sum_{j=1}^{\infty}
[(\mathbf{\rho},\mathbf{q}_{j})(\mathbf{w},\mathbf{q}_{j})
] =  \sum_{j=1}^{\infty}
[(\mathbf{\rho},\mathbf{q}_{j}) /  \nu^{2k}_j ] [ \nu^{2k}_j (\mathbf{w},\mathbf{q}_{j})
] \leq 
M_{2k}\,\| \mathcal{N}_d^k\,\mathbf{w}\|,    \eqno{(2.12)} $$ где
\begin{equation*} \label
{karo 1}
M_{2k}= \{\sum_{j=1}^{\infty}{\nu}_{j}^{-4k}(\mathbf{\rho},
\mathbf{q}_{j})^2\}^{1/2}.  
\end{equation*}

Знак равенства при заданных $(\mathbf{\rho},\mathbf{q}_{j})$
достижим. 	Значит, имеет место
\begin{lemma}
   Условие
\[M_{2k}^2= \sum_{j=1}^{\infty}{\nu}_{j}^{-4k}(\mathbf{\rho},
\mathbf{q}_{j})^2  <\infty.  \eqno(2.13)\]
		необходимо и достаточно для принадлежности
	${\rho}(\mathbf{x})$ к $( \mathcal{A}^{2k}_0)^*$.
\end{lemma}

 Величина $M_{2k}$  есть норма функционала $\rho $ в
$( \mathcal{A}^{2k}_0)^*$, которая совпадает с нормой элемента
\[ \mathcal{N}_d^{-k}\,\mathbf{f}=	(-1)^k
\sum_{j=1}^{\infty}\nu^{-2k}_j
(\mathbf{f},\mathbf{q}_{j})\mathbf{q}_{j}(\mathbf{x}) \quad \text{при} \quad  \mathbf{f} \in  \mathbf{A}^{-2k}. \eqno{(2.14)}\]

Шаг 3-й: {\it   пространство $ (\mathbf{A}^{2k}_0)^*$  отождествим с пространством $ \mathbf{A}^{-2k}$,} 

Скалярное произведение в нем	определим  как 
\[\{\mathbf{u},\mathbf{w}\}_{-2k}=
( \mathcal{N}_d^{-k}\,\mathbf{u}, \mathcal{N}_d^{-k}\,\mathbf{w}), \eqno{(2.15)}\]      а
Лемму 2.1 переформулируем 
так

\begin{theorem}
		При заданном $\mathbf{v}\in \mathcal{A}^*$ и $k\geq 1$
	уравнение $(\nabla \mathrm{div})^{2k}\,\mathbf{u}=\mathbf{v}$  разрешимо в
	пространстве  $ \mathbf{A}^{2k}$  тогда и только тогда, когда
	$\mathbf{v}\in\mathbf{A}^{-2k}$.\quad
	Его решение   $\mathbf{u}= \mathcal{N}_d^{-2k}\mathbf{v}$
	в  фактор-пространстве	$\mathcal{A}/\mathcal{A}_H$   определяется однозначно.
\end{theorem}

Действительно, если   функционал 
$\mathbf{v}\in ( \mathbf{A}^{2k}_0)^*$, то его норма  $M_{2k}<\infty$
и он принадлежит $\mathbf{A}^{-2k}$, так как 	$\{\mathbf{v},\mathbf{v}\}_{-2k}=
( \mathcal{N}_d^{-k}\,\mathbf{v},  \mathcal{N}_d^{-k}\,\mathbf{v})=M_{2k}^2$.

Ряд $ \mathcal{N}_d^{k}\,\mathbf{u}= \mathcal{N}_d^{k}\,[  \mathcal{N}_d^{-2k},\mathbf{v}]=  \mathcal{N}_d^{-k}\,\mathbf{v}$ сходится в $\mathcal{A}_{\gamma}$, так как 
$( \mathcal{N}_d^{-k}\,\mathbf{v},  \mathcal{N}_d^{-k}\,\mathbf{v})=M_{2k}^2$.
Элемент $\mathbf{u}$ принадлежит  $\mathbf{A}^{2k}$ и удовлеторяет уравнению \newline  $(\nabla \mathrm{div})^{2k}\,\mathbf{u}= \mathcal{N}_d^{2k}\,[ \mathcal{N}_d^{-2k}\,\mathbf{v}]= \mathbf{v}$,   так как квадрат его нормы
\[\{\mathbf{u},\mathbf{u}\}_{2k}=
( \mathcal{N}_d^{k}\,\mathbf{u},  \mathcal{N}_d^{k}\,\mathbf{u})=( \mathcal{N}_d^{-k}\,\mathbf{v},  \mathcal{N}_d^{-k}\,\mathbf{v})=\{\mathbf{v},\mathbf{v}\}_{-m}=M_{2k}^2<\infty. \]
Однозначноcть решения вытекает из определения и обратимости операторов   $ \mathcal{N}_d$. Теорема доказана.

Эта теорема показывает, что имеется соответсвие между пространствами 	 $ \mathbf{A}^{2k}$ и  $ \mathbf{A}^{-2k}$.
Такое же соответствие имеется между пространствами 	 $ \mathbf{W}^{m}$ и  $ \mathbf{W}^{-m}$.

{\bf Благодарности: } Академику РАН, профессору В.П.Маслову,  профессору, доктору ф.-м. наук М.Д.Рамазанову, 
профессору,	доктору ф.-м. наук С.Ю.Доброхотову, 	доктору ф.-м. наук 
Б.И.Сулейманову и  кандидату ф.-м. наук Р.Н.Гарифуллину за  поддержку.
\end{fulltext}

***


Реферат:

Операторы вихрь и градиент дивергенции в пространствах Соболева 

 Р.С.Сакс
 
{\small   Изучаются свойства  операторов вихрь    и градиент дивергенции (	$ \text{rot}$ и	$\nabla \text{div}$)
	в   пространстве $\mathbf{L}_{2}(G)$ 	
	в ограниченной области   $ G\subset  \textrm{R}^3$ с гладкой границей $\Gamma$  и в \newline   пространствах Соболева:    $ \mathbf{C}(2k, m)\equiv \mathbf{A}^{2k}(G) \oplus \mathbf{W}^m(G)$.  
	
	Пространство 	$\mathbf{L}_{2}(G) $  разлагается на ортогональные подпространства: 	
	классы  $\mathcal{A}$  и $\mathcal{B}$: 	$\mathbf{L}_{2}(G)=\mathcal{A}\oplus \mathcal{B}$. 	
	В свою очередь   
	$ \mathcal{A}= \mathcal{A}_H\oplus \mathbf{A}^0$ и $ \mathcal{B}=\mathcal{B}_H \oplus \mathbf{V}^0$,
	где 	$\mathcal{A}_H $ и $\mathcal{B}_H $- нуль-пространства операторов 	$\nabla \text{div}$   в $\mathcal{A}$  и  $ \text{rot}$  в  $\mathcal{B}$; 	
	размерности   $\mathcal{A}_H $ и $\mathcal{B}_H $ конечны и определяется топологией границы.

	Собственные поля	
	оператора	$\nabla \text{div}$  (соотв. 	$ \text{rot}$)  с ненулевыми  собственными значениями  
	используются при построении  ортонормированного базиса  в  $ \mathbf{A}^0$  (соотв. в $\mathbf{V}^0$). 	 

	Операторы   $\nabla\mathrm{div}$ и  $\mathrm{rot}$ аннулируют друг друга	  и 	проектируют $\mathbf{L}_{2}(G) $   на   $\mathcal{A}$  и $\mathcal{B}$,  причем $\mathrm{rot}\, \mathbf {u}=0$ при $\mathbf{u}\in	\mathcal{A}$,  а $\nabla\mathrm div \mathbf {v}=0$ при $\mathbf{v}\in	\mathcal{B}$  \cite{hw}.	
	
	Оператор Лапласа в  $\mathbf{L}_{2}(G)$
	выражается через них: 
	$\mathrm{\Delta} \mathbf {v}
	\equiv \nabla \mathrm{div}\,\mathbf {v}
	-(\mathrm{rot})^2\, \mathbf {v}$.	
	Поэтому 
	$\Delta^m \mathbf{v}\equiv  (\nabla \mathrm{div})^m\mathbf{v}$,  если    $  \mathbf{v}\in \mathcal{A}$,  и  $ \Delta^m \mathbf{u} \equiv(-1)^m(\mathrm{rot})^{2m}\, \mathbf{u}$,   если       $\mathbf{u}\in \mathcal{B}$, значит  $\Delta^m$-это    $ (\nabla \mathrm{div})^m$  в   $\mathcal{A}$ и  $(-1)^m(\mathrm{rot})^{2m}$- в $\mathcal{B}$ при   $ m\geq  1$.  
	
	С.Л.Соболев
	изучил  краевые задачи для  полигармонического уравнения $\Delta^m\,u=\rho$ в пространствах ${W} _2^m(\Omega)$ c обобщённой правой частью   и заложил фундамент  теории этих  пространств (см. \cite{sob},   $\S 9$ гл. 12). Его построения имеют матричные аналоги.	Вот некоторые из них.
	
	Аналоги 	 пространств  	${W}_2^{(m)}(G)$ в  классах  $\mathcal{A}$  и $\mathcal{B}$  -	это пространства   $\mathbf{A}^{2k}(G)$ и 	$\mathbf{W}^m(G)$,   порядков   $2k>0$  и $ m>0$,  а      $\mathbf{A}^{-2k}(G)$ и 	$\mathbf{W}^{-m}(G)$-  их сопряжённые пространства.
	Они определяются так:
	\[ \mathbf{A}^{2k}\equiv \{\mathbf{f}\in \mathcal{A}^0_{\gamma }, ...,(\nabla  \mathrm{div})^k \mathbf{f}\in \mathcal{A}^0_{\gamma }\}\quad \text{и} \quad  \mathbf{W}^m \equiv \{\mathbf{g}\in \mathbf{V}^0,..., (\mathrm{rot})^m \mathbf{g} \in \mathbf{V}^0\}\]
   и	образуют две шкалы (цепи) вложенных пространств:
	\[\subset \mathbf{A}^{2k}\subset...\subset\mathbf{A}^{2}\subset  \mathbf{A}^0\subset \mathbf{A}^{-2}\subset...\subset \mathbf{A}^{-2k}\subset \eqno{(0.1)}\]
	\[\subset \mathbf{W}^{m}\subset...\subset \mathbf{W}^1\subset \mathbf{V}^0\subset  \mathbf{W}^{-1}\subset...\subset \mathbf{W}^{-m}\subset     \eqno{(0.2)}\]

	В них   действуют  операторы  $\mathcal{N}_d$ и $S$ - самосопряженные расширения   операторов	$\nabla \text{div}$ и  $ \text{rot}$  в   пространства  $ \mathbf{A}^0$ и  $\mathbf{V}^{0}$.
	
	$\mathcal{N}_d$ и $S$ отображают пространства   $ \mathbf{A}^{2k}$ на  $ \mathbf{A}^{2(k-1)}$  и  $ \mathbf{W}^{m}$ на  $ \mathbf{W}^{m-1}$, соответственно,  а   операторы    $\mathcal{N}^{-1}_d$  и $S^{-1}$ -- обратно.
	
	Изучены  также отображения 	$\mathcal{N}_d^{2k}:\mathbf{A}^{2k}(G)\to \mathbf{A}^{-2k}(G)$,\quad 	$S^{2m}: \mathbf{W}^m\to \mathbf{W}^{-m}$  для любого $m\geq 1$
	и	  $k\geq 1$ и обратные  отображения 
	$\mathcal{N}_d^{-2k}$,	$S^{-2m}$   .  
	
	Доказано, что
	{\it  	уравнение $(\nabla \mathrm{div})^{2k}\,\mathbf{u}=\mathbf{v}$ 
		при заданном $\mathbf{v}$  в 	объединении  $ \mathbf{A}^{-2n}$ и $k\geq 1$ разрешимо в
		пространстве $\mathbf{A}^{2k}$ тогда и только тогда, когда
		$\mathbf{v}\in\mathbf{A}^{-2k}$.\quad
		Его решение   $\mathbf{u}= \mathcal{N}_d^{-2k}\mathbf{v}$
		в  фактор-пространстве	$\mathcal{A}/\mathcal{A}_H$   определяется однозначно. }

	Аналогично,	{\it  	при заданном $\mathbf{v}$  в 
		объединении  $\mathbf{W}^{-m}$
		уравнение $\mathrm{rot}^{2m}\,\mathbf{u}=\mathbf{v}$  разрешимо в
		пространстве $\mathbf{W}^m(G)$ тогда и только тогда, когда
		$\mathbf{v}\in \mathbf{W}^{-m}(G)$.\,
		Его решение $\mathbf{u}=S^{-2m}\mathbf{v}$ в  классе смежности   
		$ \mathcal{B}(G)/\mathcal{B}_H(G)$  определяется 
		однозначно.}

	Пары пространств из цепочек (0.1) и (02) образуют сеть пространств Соболева,  её элементы - классы  $ \mathbf{C}(2k, m)(G)\equiv \mathbf{A}^{2k}(G) \oplus \mathbf{W}^m(G)$;
	класс  $  \mathbf{C}(2k, 2k)$ совпадает с пространством Соболева    $\mathbf{H}^{2k}(G)$. 
	Они принадлежат     $\mathbf{L}_{2}(G)$, если  $k\geq 0$ и
	$m\geq 0$.
	
	Открылось широкое поле задач: изучение операторов  $(\mathrm{rot})^p$, $ (\nabla \, \mathrm{div})^p$  при  $p=1,2, ...$ и других в  сети пространств Соболева.

*

\begin{keywords}  пространство Лебега и пространства  Соболева, операторы градиент, дивергенция,  ротор, потенциальные и вихревые поля,  поля Бельтрами,	эллиптические  краевые и спектральные задачи.\end{keywords}

*
	  
Abstract: Saks Romen Semenovich

The Operators: Vortex  and  the Gradient of Divergence  in Sobolev Spaces

The properties of the vortex and the  gradient of divergence  operators 
( $ \text{rot}$	 and  $\nabla \text{div}$  ) are studied
in the space $ \mathbf {L} _{2} (G) $ 
in a bounded domain $ G \subset \textrm {R}^3 $ with a smooth boundary
 $ \Gamma$ and in the Sobolev spaces:  $ \mathbf{C}(2k, m)(G)\equiv \mathbf{A}^{2k}(G) \oplus \mathbf{W}^m(G)$. 
 
The space $ \mathbf {L}_{2} (G) $ is decomposed into orthogonal subspaces:  classes   $ \mathcal{A} $ and $ \mathcal {B} $, 
	$\mathbf{L}_{2}(G)=\mathcal{A}\oplus \mathcal{B}$.	
	In turn,  $ \mathcal{A}= \mathcal{A}_H\oplus \mathbf{A}^0$ and  $\mathcal{B}=\mathcal{B}_H \oplus \mathbf{V}^0$,
where  $\mathcal{A}_H $ и $\mathcal{B}_H $ are null spaces of operators
  	$\nabla \text{div}$ and  $ \text{rot}$  in $\mathcal{A}$  and $\mathcal{B}$;   the dimensions of   $\mathcal{A}_H $ and $\mathcal{B}_H $   are finite and  determined by the topology of the boundary.
  	
  	 In the class  $ \mathbf{A}^0$ (resp.,  In $\mathbf{V}^0$ )  an orthonormal basis is constructed from the eigenfields
  	$\mathbf{q}_{j}(\mathbf{x})$  of	$\nabla \text{div}$ operatoг (resp.,	$\mathbf{q}^{\pm }_{j}(\mathbf{x})$ of	$ \text{rot}$  operatoг) with nonzero eigenvalues    $\mu_{j}$   (resp., $\pm \lambda_{j}$ ). 	 
   Their completeness has been proven 	 \cite{ yogi, ds18}.	
  
  The operators $\nabla\mathrm{div}$ and $\mathrm{rot}$ cancel each other out and project $\mathbf{L}_{2}(G) $  onto $ \mathcal {A} $ and $ \mathcal { B} $, and $ \mathrm {rot} \, \mathbf {u} = 0 $ for $ \mathbf {u} \in \mathcal {A} $, and $ \nabla \mathrm div \mathbf {v} = 0 $ for $ \mathbf {v} \in \mathcal {B} $ \cite{hw}.
  	
Laplace matrix operator
expressed through them:
 $\mathrm{\Delta} \mathbf {v}
 \equiv \nabla \mathrm{div}\,\mathbf {v}
 -(\mathrm{rot})^2\, \mathbf {v}$.	
 Therefore
 $\Delta^m \mathbf{v}\equiv  (\nabla \mathrm{div})^m\mathbf{v}$,  if    $  \mathbf{v}\in \mathcal{A}$,  and  $ \Delta^m \mathbf{u} \equiv(-1)^m(\mathrm{rot})^{2m}\, \mathbf{u}$,   if     $\mathbf{u}\in \mathcal{B}$;   it means  $\Delta^m$ is   $ (\nabla \mathrm{div})^m$  in   $\mathcal{A}$ and it is  $(-1)^m(\mathrm{rot})^{2m}$ in $\mathcal{B}$ for all  $ m\geq  1$.  
  
S.L. Sobolev studied boundary value problems for the scalar polyharmonic equation $\Delta^m\,u=\rho$ in the spaces ${W} _2^m(\Omega)$  with a generalized right-hand side and laid the foundation for the theory of these spaces
 (\cite{sob},   $\S 9$ гл. 12).  Its constructions have matrix analogs, here are some of them.
   Analogues of the spaces ${W}_2^{(m)}(G)$ in the classes $ \mathcal {A} $ and $ \mathcal {B} $ are the space  $\mathbf{A}^{2k}(G)$   and 	$\mathbf{W}^m(G)$ of orders $ 2k> 0 $ and $ m> 0 $, and
  $ \mathbf {A}^{-2k} (G) $ and  their dual spaces $ \mathbf{W}^{- m}(G) $.
 They are defined like this:
   \[ \mathbf{A}^{2k}\equiv \{\mathbf{f}\in \mathcal{A}^0_{\gamma }, ...,(\nabla  \mathrm{div})^k \mathbf{f}\in \mathcal{A}^0_{\gamma }\}\quad and \quad  \mathbf{W}^m \equiv \{\mathbf{g}\in \mathbf{V}^0,..., (\mathrm{rot})^m \mathbf{g} \in \mathbf{V}^0\}.\]
  and form two scales (chains) of nested spaces::
\[\subset \mathbf{A}^{2k}\subset...\subset\mathbf{A}^{2}\subset  \mathbf{A}^0\subset \mathbf{A}^{-2}\subset...\subset \mathbf{A}^{-2k}\subset \eqno{(0.1)}\]
\[\subset \mathbf{W}^{m}\subset...\subset \mathbf{W}^1\subset \mathbf{V}^0\subset  \mathbf{W}^{-1}\subset...\subset \mathbf{W}^{-m}\subset     \eqno{(0.2)}\]	  

	The operators   $\mathcal{N}_d$ and  $S$  act in them; they are   self-adjoint extensions of the operators	$\nabla \text{div}$  and   $ \text{rot}$  in the spaces $ \mathbf{A}^0$  and   $\mathbf{V}^{0}$. 
	
	$\mathcal{N}_d$ and $S$	
 map the space  $ \mathbf{A}^{2k}$  to   $ \mathbf{A}^{2(k-1)}$ and  $ \mathbf{W}^{m}$  to  $ \mathbf{W}^{m-1}$,  respectively, and the operators $ \mathcal {N} ^{- 1} _d $ and $ S ^{- 1} $ - back.

Examined mappings
$ S^{2m}: \mathbf {W}^m \to \mathbf {W}^{- m} $ and
$ S^{- 2m}: \mathbf {W}^{- m} \to \mathbf {W}^m $ for any $ m \geq 1 $.
It is proved that
{\it for a given $ \mathbf{v} $ in 
	union of  	 $\mathbf{W}^{-m}$
	 and $ m \geq 1 $
	the equation $ \mathrm {rot}^{2m} \, \mathbf {u} =\mathbf {v} $ is solvable in
	space $ \mathbf {W} ^m (G) $ if and only if
	$ \mathbf {v} \in \mathbf {W}^{- m} (G) $. \,
	Its solution $ \mathbf{u} = S^{- 2m} \mathbf {v} $ in the coset  $ \mathcal {B} (G)/ \mathcal {B}_H$  is  uniquely defined}.

Similar properties of mappings  $\mathcal{N}_d^{2k}:\mathbf{A}^{2k}(G)\to \mathbf{A}^{-2k}(G)$  and $\mathcal{N}_d^{-2k}$  
 for $ k \geq 1 $ are proved in Theorem 4.

Pairs of  spaces from chains (0.1) and (02) form a net of Sobolev spaces, its elements are classes  $ \mathbf{C}(2k, m)(G)\equiv \mathbf{A}^{2k}(G) \oplus \mathbf{W}^m(G)$;
the class $  \mathbf{C}(2k, 2k)$coincides with the Sobolev space   $\mathbf{H}^{2k}(G)$. 
They belong to     $\mathbf{L}_{2}(G)$,  if  $k\geq 0$  and $m\geq 0$.

A wide field of problems has opened up: studying the operators $(\mathrm{rot})^p$, $ (\nabla \, \mathrm{div})^p$  for $ p = 1,2, ...,$  and others in the network Sobolev spaces.

*

{\bf Keywords:} Lebesgue space and Sobolev spaces, operators: gradient, divergence, curle, potential and vortex fields, Beltrami fields, elliptic boundary value and spectral problems.

*

Сакс Ромэн Семенович

старший научный сотрудник

 Институт Математики с ВЦ УФИЦ РАН 

450077, г. Уфа, ул. Чернышевского, д.112 

телефон: (347)272-59-36 \quad  (347) 273-34-12

факс: (347) 272-59-36

 телефон дом.: (347) 273-84-69
 
  моб.+7 917 379 75 38

 e-mail: romen-saks@yandex.ru

\end{document}